\documentclass[11pt]{article}
	
	\newcommand{\mainTitle}{Congruence identities of regularized multiple zeta values involving a pair of index sets}

	\newcommand{\authorName}{MACHIDE, Tomoya}
	\newcommand{\organizationNameFst}{National Institute of Informatics}
	\newcommand{\placeAddressFst}{2-1-2 Hitotsubashi, Chiyoda-ku, Tokyo 101-8430, Japan}
	\newcommand{\emailAddressFst}{machide@nii.ac.jp}
	\newcommand{\organizationNameScd}{JST, ERATO, Kawarabayashi Large Graph Project}
	\newcommand{\departmentNameScd}{c/o Global Research Center for Big Data Mathematics}
	\newcommand{\placeAddressScd}{NII, 2-1-2 Hitotsubashi, Chiyoda-ku, Tokyo 101-8430, Japan}
	
\usepackage{graphicx}
\usepackage{geometry}               
\usepackage{ascmac}
\usepackage{amsmath}
\usepackage{amssymb}
\usepackage{amsthm}
\usepackage[hypertex, colorlinks=true, bookmarks=true]{hyperref}
\usepackage{arydshln}

	\DeclareMathOperator*{\OPlus}{\bigoplus}

	\newcommand{\nbk}[3]{#1#3#2}		
	\newcommand{\bgbk}[3]{\bigl{#1}#3\bigr{#2}}	
	\newcommand{\Bgbk}[3]{\Bigl{#1}#3\Bigr{#2}}			
	\newcommand{\bggbk}[3]{\biggl{#1}#3\biggr{#2}}			
	\newcommand{\Bggbk}[3]{\Biggl{#1}#3\Biggr{#2}}
	\newcommand{\autobk}[3]{\left#1#3\right#2}
	\newcommand{\nbkD}[5]{#1#2#5#3#4}		
	\newcommand{\bgbkD}[5]{\bigl{#1}\bigl{#2}#5\bigr{#3}\bigr{#4}}	
	\newcommand{\BgbkD}[5]{\Bigl{#1}\Bigl{#2}#5\Bigr{#3}\Bigr{#4}}	
	\newcommand{\bggbkD}[5]{\biggl{#1}\biggl{#2}#5\biggr{#3}\biggr{#4}}	
	\newcommand{\BggbkD}[5]{\Biggl{#1}\Biggl{#2}#5\Biggr{#3}\Biggr{#4}}	
	\newcommand{\autobkD}[5]{\left#1\left#2#5\right#3\right#4}	
	\newcommand{\mcbk}[4][?]{\ifx n#1\nbk{#2}{#3}{#4}\else\ifx b#1\bgbk{#2}{#3}{#4}\else\ifx B#1\Bgbk{#2}{#3}{#4}\else\ifx g#1\bggbk{#2}{#3}{#4}\else\ifx G#1\Bggbk{#2}{#3}{#4}\else\ifx a#1\autobk{#2}{#3}{#4}\else\ifx !#1{#4}\else#4\fi\fi\fi\fi\fi\fi\fi}
	\newcommand{\mcbkD}[4][?]{\ifx n#1\nbkD{#2}{#2}{#3}{#3}{#4}\else\ifx b#1\bgbkD{#2}{#2}{#3}{#3}{#4}\else\ifx B#1\BgbkD{#2}{#2}{#3}{#3}{#4}\else\ifx g#1\bggbkD{#2}{#2}{#3}{#3}{#4}\else\ifx G#1\BggbkD{#2}{#2}{#3}{#3}{#4}\else\ifx a#1\autobkD{#2}{#2}{#3}{#3}{#4}\else\ifx !#1{#4}\else#4\fi\fi\fi\fi\fi\fi\fi}
	\newcommand{\nsgsb}[1]{#1}		
	\newcommand{\bgsgsb}[1]{\big{#1}}	
	\newcommand{\Bgsgsb}[1]{\Big{#1}}			
	\newcommand{\bggsgsb}[1]{\bigg{#1}}			
	\newcommand{\Bggsgsb}[1]{\Bigg{#1}}
	\newcommand{\mcsgsb}[2][?]{\ifx n#1\nsgsb{#2}\else\ifx b#1\bgsgsb{#2}\else\ifx B#1\Bgsgsb{#2}\else\ifx g#1\bggsgsb{#2}\else\ifx G#1\Bggsgsb{#2}\else#2\fi\fi\fi\fi\fi}
	\newcommand{\myEqSpace}{\,} 	\newlength{\myEqSpaceLen} 	
	\newcommand{\mLt}[1]{\widetilde{#1}}

	\newcommand{\bkR}[2][a]{\mcbk[#1]{(}{)}{#2}}						
	\newcommand{\bkS}[2][a]{\mcbk[#1]{[}{]}{#2}}						
	\newcommand{\bkB}[2][a]{\mcbk[#1]{\{}{\}}{#2}}						
	\newcommand{\bkA}[2][a]{\mcbk[#1]{\langle}{\rangle}{#2}}				
		
	\newcommand{\nFc}[3][n]{#2\bkR[#1]{#3}}
	
	\newcommand{\idFc}[4][n]{\id{#2}{#3}\bkR[#1]{#4}}
	\newcommand{\pwFc}[4][n]{\pw{#2}{#3}\bkR[#1]{#4}}
	\newcommand{\ipFc}[5][n]{\ip{#2}{#3}{#4}\bkR[#1]{#5}}
	
		\newcommand{\Fc}{\nFc}
			
	\newcommand{\gam}{\gamma}

	\newcommand{\ep}{\varepsilon}
	
	\newcommand{\sig}{\sigma}

	\newcommand{\vphi}{\varphi}
	\newcommand{\modd}[1][\ ]{\mathrm{mod}#1}
	\newcommand{\mo}{(-1)}

	\newcommand{\SetO}[2][n]{\bkB[#1]{#2}}
	
		\newcommand{\Set}{\SetO}

	\newcommand{\setZ}{\mathbb{Z}} 
	\newcommand{\setQ}{\mathbb{Q}}	
		
	\newcommand{\setR}{\mathbb{R}}

	\newcommand{\gpSym}[2][?]{S^{(#2)}}
	
	\newcommand{\gpKleinF}[1][?]{V}
		
	\newcommand{\gpGL}[2]{GL_{#1}(#2)}
	
	\newcommand{\gpu}[1][?]{\ifx?#1e\else e_{#1}\fi}

	\newcommand{\vPack}[1][10]{\vspace{-#1pt}}
	
	\newcommand{\lnA}[1][]{&  &}
	\newcommand{\lnP}[1]{\myEqSpace#1\myEqSpace}
	\newcommand{\lnAP}[2][]{& #2 &}
	\newcommand{\lnAH}[1][\nonumber]{#1\\ & &}

		\newcommand{\slnAH}[1][?]{\\}
		
	\newcommand{\refEq}[1]{(\ref{#1})}	

	\newcommand{\pcstSpForRefThm}{\;}		\newcommand{\refHL}[2]{#1\pcstSpForRefThm\ref{#2}}		\newcommand{\refHLm}[2]{#1\pcstSpForRefThm#2}
	\newcommand{\refThm}[2][?]{\ifx?#1\refHL{Theorem}{#2}\else\ifx s#1\refHL{Theorems}{#2}\else{[argument error]}\fi\fi}
	\newcommand{\refProp}[2][?]{\ifx?#1\refHL{Proposition}{#2}\else\ifx s#1\refHL{Propositions}{#2}\else{[argument error]}\fi\fi}
	\newcommand{\refLem}[2][?]{\ifx?#1\refHL{Lemma}{#2}\else\ifx s#1\refHL{Lemmas}{#2}\else{[argument error]}\fi\fi}
	\newcommand{\refCor}[2][?]{\ifx?#1\refHL{Corollary}{#2}\else\ifx s#1\refHL{Corollaries}{#2}\else{[argument error]}\fi\fi}
	\newcommand{\refDef}[2][?]{\ifx?#1\refHL{Definition}{#2}\else\ifx s#1\refHL{Definitions}{#2}\else{[argument error]}\fi\fi}
	\newcommand{\refRem}[2][?]{\ifx?#1\refHL{Remark}{#2}\else\ifx s#1\refHL{Remarks}{#2}\else{[argument error]}\fi\fi}
	\newcommand{\refTab}[2][?]{\ifx?#1\refHL{Table}{#2}\else\ifx s#1\refHL{Tables}{#2}\else{[argument error]}\fi\fi}
	\newcommand{\refSec}[2][?]{\ifx?#1\refHL{Section}{#2}\else\ifx s#1\refHL{Sections}{#2}\else{[argument error]}\fi\fi}
		\newcommand{\refSect}{\refSec}
	\newcommand{\refThmA}[1]{\refHLm{Theorem}{#1}}
	\newcommand{\refPropA}[1]{\refHLm{Proposition}{#1}}
	\newcommand{\refLemA}[1]{\refHLm{Lemma}{#1}}

	\newcommand{\refSecA}[1]{\refHLm{Section}{#1}}
		\newcommand{\refSectA}{\refSecA}
		
	\newcommand{\Vc}[2][?]{\ifx ?#1\vec{#2}\else\ifx l#1\overrightarrow{#2}\else\ifx b#1{\bf#2}\else[error]\fi\fi\fi}
	
	\newcommand{\nopF}[3][?]{\ifx s#1#2/#3\else\ifx b#1(#2)/(#3)\else\ifx d#1\dfrac{#2}{#3}\else\ifx t#1\tfrac{#2}{#3}\else\frac{#2}{#3}\fi\fi\fi\fi}

		\newcommand{\opF}{\nopF}

	\newcommand{\mopI}[2][?]{\bkR[#1]{#2}^{-1}}
	
			\newcommand{\mopTPm}[2][?]{\mbox{}^t\bkR[#1]{#2}}	\newcommand{\mopTP}{\mopTPm}
	\newcommand{\myVopLetter}{$\cdot$} \newlength{\myVopLetterHeight}   

	\newcommand{\racFT}[3][?]{#2\vert\bkR[#1]{#3}}						\newcommand{\racFrrT}[3][?]{(#2\vert\bkR[#1]{#3})}	
												\newcommand{\racF}{\racFT}			\newcommand{\racFrr}{\racFrrT}
	\newcommand{\racFATh}[4][n]{\bkR[#1]{#2\vert#3}(#4)}				
												\newcommand{\racFA}{\racFATh}

	\newcommand{\pw}[3][?]{\ifx!#3{#2}^{#3}\else#2^{#3}\fi}
	\newcommand{\id}[3][?]{#2_{#3}}
	\newcommand{\ip}[4][?]{{#2}_{#3}^{#4}}
	\newcommand{\pwR}[3][a]{\ifx!#1{\bkR[#1]{#2}}^{#3}\else\bkR[#1]{#2}^{#3}\fi}
	\newcommand{\pwB}[3][a]{\ifx!#1{\bkB[#1]{#2}}^{#3}\else\bkB[#1]{#2}^{#3}\fi}
	\newcommand{\pwS}[3][a]{\ifx!#1{\bkS[#1]{#2}}^{#3}\else\bkS[#1]{#2}^{#3}\fi}
					
	\newcommand{\mpIlett}{id}		\newcommand{\mpI}[1][?]{\ifx?#1\mpIlett\else \mpIlett_{#1}\fi}	
	
	\newcommand{\tpT}[3][a]{ {#2}\atop \bkR[#1]{#3} }

	\newcommand{\nSmO}[2][?]{\ifx l#1\sum\limits_{#2}\else\ifx t#1{\textstyle\sum\limits_{#2}}\else\sum_{#2}\fi\fi}
	\newcommand{\nSmT}[3][?]{\ifx l#1\sum\limits_{#2}^{#3}\else\if t#1{\textstyle\sum\limits_{#2}^{#3}}\else\sum_{#2}^{#3}\fi\fi}	
	\newcommand{\nSmN}[1][?]{\ifx l#1\sum\limits\else\ifx t#1{\textstyle\sum\limits}\else\sum\fi\fi}
	\newcommand{\pSm}[2][?]{\ifx t#1 \sum_{#2}^{\prime} \else \sideset{}{^\prime}\sum_{#2} \fi}
	\newcommand{\pSmT}[3][?]{\ifx t#1 \sum_{#2}^{\prime#3} \else \sideset{}{^\prime}\sum_{#2}^{#3} \fi}	
	\newcommand{\pSmN}[1][?]{\ifx t#1 \sum^{\prime} \else \sideset{}{^\prime}\sum \fi}
	\newcommand{\dSm}[2][?]{\ifx t#1 \sum_{#2}^{\dagger} \else \sideset{}{^\dagger}\sum_{#2} \fi}
	\newcommand{\dSmT}[3][?]{\ifx t#1 \sum_{#2}^{\dagger#3} \else \sideset{}{^\dagger}\sum_{#2}^{#3} \fi}	
	\newcommand{\dSmN}[1][?]{\ifx t#1 \sum^{\dagger} \else \sideset{}{^\dagger}\sum \fi}
	\newcommand{\tpTSm}[3][?]{\nSmO[#1]{\tpT{#2}{#3}}}

		\newcommand{\Sm}{\nSmO}			\newcommand{\SmT}{\nSmT}			\newcommand{\SmN}{\nSmN}
		\newcommand{\tpSm}{\tpTSm}
		
	\newcommand{\nPd}[2][?]{\ifx l#1 \prod\limits_{#2} \else \prod_{#2} \fi}
	\newcommand{\nPdT}[3][?]{\ifx l#1 \prod\limits_{#2}^{#3} \else \prod_{#2}^{#3} \fi}

	\newcommand{\nOPs}[2][?]{\ifx l#1 \OPlus\limits_{#2} \else \OPlus_{#2} \fi}
	\newcommand{\nOPsT}[3][?]{\ifx l#1 \OPlus\limits_{#2}^{#3} \else \OPlus_{#2}^{#3} \fi}	
	\newcommand{\pOPs}[2][?]{\ifx t#1 \OPlus_{#2}^{\prime} \else \sideset{}{^\prime}\OPlus_{#2} \fi}
	\newcommand{\pOPsT}[3][?]{\ifx t#1 \OPlus_{#2}^{\prime#3} \else \sideset{}{^\prime}\OPlus_{#2}^{#3} \fi}	
	\newcommand{\tpTOPs}[3][?]{ \nOPs[#1]{ \tpT{#2}{#3} } }

		\newcommand{\tpOPs}{\tpTOPs}
		
	\newcommand{\nIs}[2][?]{\ifx l#1 \bigcap\limits_{#2}\else\ifx b#1 \bigcap_{#2}\else{\textstyle\bigcap\limits_{#2}}\fi\fi}
	\newcommand{\nIsT}[3][?]{\ifx l#1 \bigcap\limits_{#2}^{#3}\else\ifx b#1 \bigcap_{#2}^{#3}\else{\textstyle\bigcap\limits_{#2}^{#3}}\fi\fi}	
	\newcommand{\pIs}[2][?]{\ifx t#1 \bigcap_{#2}^{\prime} \else \sideset{}{^\prime}\bigcap_{#2} \fi}
	\newcommand{\pIsT}[3][?]{\ifx t#1 \bigcap_{#2}^{\prime#3} \else \sideset{}{^\prime}\bigcap_{#2}^{#3} \fi}

	\newcommand{\nUn}[2][?]{\ifx l#1 \bigcup\limits_{#2}\else\ifx b#1 \bigcup_{#2}\else{\textstyle\bigcup\limits_{#2}}\fi\fi}
	\newcommand{\nUnT}[3][?]{\ifx l#1 \bigcup\limits_{#2}^{#3}\else\ifx b#1 \bigcup_{#2}^{#3}\else{\textstyle\bigcup\limits_{#2}^{#3}}\fi\fi}	
	\newcommand{\pUn}[2][?]{\ifx t#1 \bigcup_{#2}^{\prime} \else \sideset{}{^\prime}\bigcup_{#2} \fi}
	\newcommand{\pUnT}[3][?]{\ifx t#1 \bigcup_{#2}^{\prime#3} \else \sideset{}{^\prime}\bigcup_{#2}^{#3} \fi}

	\newcommand{\nLm}[2][?]{\ifx l#1 \lim\limits_{#2} \else \lim_{#2} \fi}
		
	\newcommand{\glcondEnvLineHead}[1]{ \ifx*#1 \begin{eqnarray*} \else \begin{eqnarray}  \label{#1} \fi }
	\newcommand{\glcondEnvLineTail}[1]{ \ifx*#1 \end{eqnarray*} \else \end{eqnarray} \fi }
	\newcommand{\glcondDis}[1]{\ifx d#1 \displaystyle \fi}
	\newcommand{\glcmdEqShift}{\hspace{-20pt}}
	\newcommand{\glcmdHLineCWiden}{\rule{0cm}{15pt}}	\newcommand{\glcdH}{\glcmdHLineCWiden}
	\newcommand{\lccondPar}[1]{\ifx#1p \\ \fi}

		\newcommand{\envMO}[2][*]{$\ifx d#1 \displaystyle \fi#2$}
		\newcommand{\envMT}[3][*]{$\ifx d#1 \displaystyle \fi#2=#3$}
		\newcommand{\envMTDef}[3][*]{$\ifx d#1 \displaystyle \fi#2:=#3$}
		\newcommand{\envMTPt}[4][*]{$\ifx d#1 \displaystyle \fi#3#2#4$}
			\newcommand{\envM}{\envMT}
			
			\newcommand{\envMPt}{\envMTPt}
		\newcommand{\envMTh}[4][*]{$\ifx d#1 \displaystyle \fi#2=#3=#4$}
		\newcommand{\envMThPt}[5][*]{$\ifx d#1 \displaystyle \fi#3#2#4#2#5$}
		\newcommand{\envMF}[5][*]{$\ifx d#1 \displaystyle \fi#2=#3=#4=#5$}
		\newcommand{\envMFPt}[6][*]{$\ifx d#1 \displaystyle \fi#3#2#4#2#5#2#6$}
		
	\newcommand{\envMLineT}[3][*]{ \ifx*#1 \begin{multline*} #2\lnP{=}#3\end{multline*} \else \begin{multline} \label{#1} #2\lnP{=}#3\end{multline} \fi }
	\newcommand{\envMLineTDef}[3][*]{ \ifx*#1 \begin{multline*} #2\lnP{:=}#3\end{multline*} \else \begin{multline} \label{#1} #2\lnP{:=}#3\end{multline} \fi }
	\newcommand{\envMLineTPt}[4][*]{ \ifx*#1 \begin{multline*} #3\lnP{#2}#4\end{multline*} \else \begin{multline} \label{#1} #3\lnP{#2}#4\end{multline} \fi }

		\newcommand{\envHLineT}[3][*]{ \glcondEnvLineHead{#1} #2&=&#3\glcondEnvLineTail{#1} }
		\newcommand{\envHLineTDef}[3][*]{ \glcondEnvLineHead{#1} #2&:=&#3\glcondEnvLineTail{#1} }
		\newcommand{\envHLineTPt}[4][*]{\glcondEnvLineHead{#1} #3&#2&#4\glcondEnvLineTail{#1}}
		
			\newcommand{\envHLine}{\envHLineT}
			\newcommand{\envHLineDef}{\envHLineTDef}
			\newcommand{\envHLinePt}{\envHLineTPt}
			
		\newcommand{\envPLineTNme}[3][*]{\begin{eqnarray} #2\\#3\end{eqnarray}}

		\newcommand{\envOTLine}[4][*]{\glcondEnvLineHead{#1} #2\lnAP{=}#3\lnP{=}#4\glcondEnvLineTail{#1}}

	\newcommand{\lcparaCase}{\vspace{3pt}}

	\newcommand{\envSMatT}[3][a]{\autobk{(}{)}{\begin{smallmatrix}#2\\#3\end{smallmatrix}}}
	\newcommand{\envMatT}[3][a]{\autobk{(}{)}{\begin{matrix}#2\\#3\end{matrix}}}

	\newcommand{\envMatTh}[4][a]{ \autobk{(}{)}{\begin{matrix}#2\\#3\\#4\end{matrix}} }
	
	\newcommand{\envMatF}[5][a]{ \autobk{(}{)}{\begin{matrix}#2\\#3\\#4\\#5\end{matrix}} }

		\newcommand{\SMatT}{\envSMatT}

		\newcommand{\MatT}{\envMatT}
		
		\newcommand{\MatTh}{\envMatTh}
		
		\newcommand{\MatF}{\envMatF}

	\newcommand{\matu}[1][?]{\ifx#1?I\else I_{#1}\fi}

	\newcommand{\exx}[2][n]{ \Fc[#1]{\exp}{#2} }

	\newcommand{\fcZlett}{\zeta}				\newcommand{\fcZO}[2][n]{\Fc[#1]{\fcZlett}{#2}}	
																\newcommand{\fcZ}{\fcZO}
	\newcommand{\lgP}[3][n]{\idFc[#1]{Li}{#2}{#3}}
	
	\newcommand{\ctG}[1][?]{\gam}

	\newcommand{\cTxT}[2]{\textcolor{#1}{#2}}
	
		\newcommand{\cTx}{\cTxT}
	\newcommand{\alTx}[1]{\cTx{red}{#1}}
	\newcommand{\rmTx}[1]{\cTx{blue}{#1}}
	\newcommand{\ntTx}[1]{\cTx{Magenta}{#1}}
	
	\newcommand{\cmoTx}[1]{\cTx{Gray}{#1}}
	\newcommand{\sTx}[2][?]{ \ifx t#1{\tiny #2} \else \ifx s#1{\scriptsize #2} \else \ifx f#1{\footnotesize #2} \else \ifx S#1{\small #2} \else \ifx n#1{\normalsize #2} \else \ifx l#1{\large #2} \else \ifx L#1{\Large #2} \else \ifx R#1{\LARGE #2} \else \ifx h#1{\huge #2} \else \ifx H#1{\Huge #2} \else \ifx ?#1 #2 \else #2 \fi\fi\fi\fi\fi\fi\fi\fi\fi\fi\fi }
	\newcommand{\bfTx}[1]{{\bf#1}}
	\newcommand{\respTx}[1]{(resp. #1)}
	
	\newcommand{\osMTx}[3][?]{\overset{#3}{#2}}
	\newcommand{\usbMTx}[3][?]{\underset{#2}{\underbrace{#3}}}

		\newcommand{\osTx}{\osMTx}
		\newcommand{\usbTx}{\usbMTx}
		
	\newcommand{\raMTx}[3][?]{\raisebox{#2pt}[0pt][0pt]{$\ifx d#1\displaystyle\fi#3$}}
	
	\newcommand{\roMTx}[3][?]{\rotatebox[origin=c]{#2}{$#3$}}

		\newcommand{\envHLineTCl}[3][a]{ \ifx a#1\alTx{\envHLineT{#2}{#3}} \else \ifx r#1 \rmTx{\envHLineT{#2}{#3}} \else\ifx n#1 \ntTx{\envHLineT{#2}{#3}}\else\ifx c#1 \cmoTx{\envHLineT{#2}{#3}} \else \text{[argument error]} \fi\fi\fi\fi \vPack[18] }

		\newcommand{\envHLineCSClPart}[8][?]{\ifx*#1 \begin{eqnarray*} \else \begin{eqnarray}  \label{#1}  \fi \alTx{#3}&#2&\alTx{#4}\\\glcdH#5&#2&#6\nonumber\\\glcdH\alTx{#7}&#2&\alTx{#8}\nonumber\glcondEnvLineTail{*}}
		
			\newcommand{\HLineCTCl}[3][?]{\alTx{#2}&=&\alTx{#3}\nonumber \ifx#1p \\\glcdH \fi}
			\newcommand{\HLineCTClDef}[3][?]{\alTx{#2}&:=&\alTx{#3}\nonumber \ifx#1p \\\glcdH \fi}
			\newcommand{\HLineCFCl}[5][?]{\alTx{#2}&=&\alTx{#3}\nonumber\\\glcdH#4&=&#5\nonumber \ifx#1p \\\glcdH \fi}
			\newcommand{\HLineCFClDef}[5][?]{\alTx{#2}&:=&\alTx{#3}\nonumber\\\glcdH#4&:=&#5\nonumber \ifx#1p \\\glcdH \fi}
			\newcommand{\HLineCSCl}[7][?]{\alTx{#2}&=&\alTx{#3}\nonumber\\\glcdH#4&=&#5\nonumber\\\glcdH\alTx{#6}&=&\alTx{#7}\nonumber\ifx#1p \\\glcdH \fi}
			\newcommand{\HLineCSClDef}[7][?]{\alTx{#2}&:=&\alTx{#3}\nonumber\\\glcdH#4&:=&#5\nonumber\\\glcdH\alTx{#6}&:=&\alTx{#7}\nonumber\ifx#1p \\\glcdH \fi}									
			\newcommand{\HLineCECl}[9][?]{\alTx{#2}&=&\alTx{#3}\nonumber\\\glcdH#4&=&#5\nonumber\\\glcdH\alTx{#6}&=&\alTx{#7}\nonumber\\\glcdH#8&=&#9\nonumber\ifx#1p \\\glcdH \fi}
			\newcommand{\HLineCEClDef}[9][?]{\alTx{#2}&:=&\alTx{#3}\nonumber\\\glcdH#4&:=&#5\nonumber\\\glcdH\alTx{#6}&:=&\alTx{#7}\nonumber\\\glcdH#8&:=&#9\nonumber\ifx#1p \\\glcdH \fi}
		
	\newcommand{\envCenter}[2][*]{\ifx*#1\begin{center}\else\begin{center}[#1]\fi #2\end{center}}
	\newcommand{\envFlushleft}[2][*]{\ifx*#1\begin{flushleft}\else\begin{flushleft}[#1]\fi #2\end{flushleft}}
	\newcommand{\envFlushright}[2][*]{\ifx*#1\begin{flushright}\else\begin{flushright}[#1]\fi #2\end{flushright}}
		
	\newcommand{\envItemIm}[2][*]{\ifx*#1\begin{itemize}\else\begin{itemize}[#1]\fi #2\end{itemize}}
	\newcommand{\envItemDp}[2][*]{\ifx*#1\begin{description}\else\begin{description}[#1]\fi #2\end{description}}
	\newcommand{\envItemEm}[2][*]{\ifx*#1\begin{enumerate}\else\begin{enumerate}[#1]\fi #2\end{enumerate}}

	\newcommand{\envMultCol}[3][*]{\ifx1#2#3\else\begin{multicols}{#2}\ifx*#1\else\mbox{}\vspace{-#1pt}\fi#3\end{multicols}\fi}
	
\theoremstyle{plain}
\newtheorem{theorem}{THEOREM}[section]
\newtheorem{proposition}[theorem]{PROPOSITION}
\newtheorem{lemma}[theorem]{LEMMA}
\newtheorem{corollary}[theorem]{COROLLARY}
\theoremstyle{definition}

\theoremstyle{remark}
\newtheorem{remark}[theorem]{REMARK}
\theoremstyle{plain}

\theoremstyle{definition}

\theoremstyle{remark}

\theoremstyle{plain}

\theoremstyle{definition}

\theoremstyle{remark}

\allowdisplaybreaks[4]
\numberwithin{equation}{section}

	\newcommand{\lccondBibitem}[3][]{ \if ?#2 \bibitem{#3} \else \bibitem[#2]{#3} \fi}
	\newcommand{\refPaper}[8][?]{
			\lccondBibitem{#1}{#2}
				#3,			
				\emph{#4}, 	
				#5\ 			
				{\bf #6},		
				#7,			
				#8.			
		}
	\newcommand{\refPreprint}[6][?]{
			\lccondBibitem{#1}{#2}
				#3,			
				\emph{#4}, 	
				preprint; #5,			
				#6.			
		}
	\newcommand{\refPaperRep}[9][?]{
			\lccondBibitem{#1}{#2}
				#3,			
				\emph{#4}, 	
				#5\ 			
				{\bf #6},		
				#7,			
				#8			
				; reprinted in #9	
		}

	\newcommand{\refPaperAlm}[5][?]{
			\lccondBibitem{#1}{#2}
				#3,	 		
				\emph{#4}, 	
				#5		
		}

	\newcommand{\etalTx}[2][?]{#2 \emph{et al.}\!}

	\newcommand{\glcondEnvLineTailPd}[1]{.\ifx*#1 \end{eqnarray*} \else \end{eqnarray} \fi  }
	\newcommand{\glcondEnvLineTailCm}[1]{,\ifx*#1 \end{eqnarray*} \else \end{eqnarray} \fi }
	\newcommand{\prcondEnvEqSpHead}[1]{ \ifx*#1 \begin{equation*}[ERROR] \else \begin{equation}  \label{#1} \fi  }
	\newcommand{\prcondEnvEqSpTail}[1]{\ifx*#1 [ERROR]\end{equation*} \else \end{equation} \fi }

	\newcommand{\envProof}[2][?]{ \par\mbox{}\vspace{-5pt}\\ \ifx?#1\emph{Proof.}\else\emph{Proof of #1.}\fi \ #2 \hfill $\Box$\\ \par}
	
		\newcommand{\envLineThCm}[4][*]{ \glcondEnvLineHead{#1} & &\glcmdEqShift#2\nonumber\\&=&#3 \\&=&#4\nonumber \glcondEnvLineTailCm{#1} }

		\newcommand{\envLineFPd}[5][*]{ \glcondEnvLineHead{#1} & &\glcmdEqShift#2\nonumber\\&=&#3\nonumber \\&=&#4 \\&=&#5\nonumber \glcondEnvLineTailPd{#1} }
		
		\newcommand{\envHLineTPd}[3][*]{ \glcondEnvLineHead{#1} #2&=&#3\glcondEnvLineTailPd{#1} }
		\newcommand{\envHLineTDefPd}[3][*]{ \glcondEnvLineHead{#1} #2&:=&#3\glcondEnvLineTailPd{#1} }
		\newcommand{\envHLineTCm}[3][*]{ \glcondEnvLineHead{#1} #2&=&#3\glcondEnvLineTailCm{#1} }
		\newcommand{\envHLineTCmDef}[3][*]{ \glcondEnvLineHead{#1} #2&:=&#3\glcondEnvLineTailCm{#1} }
		\newcommand{\envHLineTCmPt}[4][*]{\glcondEnvLineHead{#1} #3&#2&#4\glcondEnvLineTailCm{#1}}					
		\newcommand{\envHLineTPdPt}[4][*]{\glcondEnvLineHead{#1} #3&#2&#4\glcondEnvLineTailPd{#1}}

			\newcommand{\envHLinePd}{\envHLineTPd}
			\newcommand{\envHLineDefPd}{\envHLineTDefPd}
			\newcommand{\envHLineCm}{\envHLineTCm}
			\newcommand{\envHLineCmDef}{\envHLineTCmDef}
			\newcommand{\envHLineCmPt}{\envHLineTCmPt}
			\newcommand{\envHLinePdPt}{\envHLineTPdPt}

		\newcommand{\envHLineThCmPt}[5][*]{\glcondEnvLineHead{#1} #3&#2&#4\nonumber\\&#2&#5  \glcondEnvLineTailCm{#1}}
		\newcommand{\envHLineThPdPt}[5][*]{\glcondEnvLineHead{#1} #3&#2&#4\nonumber\\&#2&#5  \glcondEnvLineTailPd{#1}}

		\newcommand{\envHLineFPdPt}[6][*]{\glcondEnvLineHead{#1} #3&#2&#4\nonumber\\&#2&#5 \\&#2&#6\nonumber\glcondEnvLineTailPd{#1}}

		\newcommand{\envHLineCFCm}[5][*]{\glcondEnvLineHead{#1} #2&=&#3,\nonumber\\\glcdH#4&=&#5\glcondEnvLineTailCm{#1}}

		\newcommand{\envHLineCFCmNme}[5][*]{\begin{eqnarray} #2&=&#3,\\\glcdH#4&=&#5 \glcondEnvLineTailCm{?} }
		\newcommand{\envHLineCFNmePd}[5][*]{\begin{eqnarray} #2&=&#3,\\\glcdH#4&=&#5 \glcondEnvLineTailPd{?} }
		\newcommand{\envHLineCFCmDefNme}[5][*]{\begin{eqnarray} #2&:=&#3,\\\glcdH#4&:=&#5 \glcondEnvLineTailCm{?} }
		\newcommand{\envHLineCFDefNmePd}[5][*]{\begin{eqnarray} #2&:=&#3,\\\glcdH#4&:=&#5 \glcondEnvLineTailPd{?} }
		\newcommand{\envHLineCFCmPt}[6][*]{\glcondEnvLineHead{#1} #3&#2&#4,\nonumber\\\glcdH#5&#2&#6\glcondEnvLineTailCm{#1}}

		\newcommand{\envHLineCFNmePdPt}[6][*]{\begin{eqnarray}#3&#2&#4,\\\glcdH#5&#2&#6\glcondEnvLineTailPd{?}}
		\newcommand{\envHLineCFCmNmePt}[6][*]{\begin{eqnarray}#3&#2&#4,\\\glcdH#5&#2&#6\glcondEnvLineTailCm{?}}
		\newcommand{\envHLineCFNmePdPte}[7][*]{\begin{eqnarray}#2&#3&#4,\\\glcdH#5&#6&#7\glcondEnvLineTailPd{?}}
		\newcommand{\envHLineCFCmNmePte}[7][*]{\begin{eqnarray}#2&#3&#4,\\\glcdH#5&#6&#7\glcondEnvLineTailCm{?}}

		\newcommand{\envHLineCFLqqPd}[5][*]{\envPLinePd[#1]{#2\lnP{=}#3,\qquad#4\lnP{=}#5}}
		\newcommand{\envHLineCFLaCm}[5][*]{\envPLineCm[#1]{#2\lnP{=}#3\quad\text{and}\quad#4\lnP{=}#5}}		
							
		\newcommand{\envHLineCFLaaCm}[5][*]{\envPLineCm[#1]{#2\lnP{=}#3\qquad\text{and}\qquad#4\lnP{=}#5}}		
		\newcommand{\envHLineCFLaaPd}[5][*]{\envPLinePd[#1]{#2\lnP{=}#3\qquad\text{and}\qquad#4\lnP{=}#5}}

		\newcommand{\envHLineCSNmePd}[7][*]{\begin{eqnarray} #2&=&#3,\\\glcdH#4&=&#5,\\\glcdH#6&=&#7\glcondEnvLineTailPd{?}}
		\newcommand{\envHLineCSDefNmePd}[7][*]{\begin{eqnarray} #2&:=&#3,\\\glcdH#4&:=&#5,\\\glcdH#6&:=&#7\glcondEnvLineTailPd{?}}
		\newcommand{\envHLineCSCmNme}[7][*]{\begin{eqnarray} #2&=&#3,\\\glcdH#4&=&#5,\\\glcdH#6&=&#7\glcondEnvLineTailCm{?}}
		\newcommand{\envHLineCSCmDefNme}[7][*]{\begin{eqnarray} #2&:=&#3,\\\glcdH#4&:=&#5,\\\glcdH#6&:=&#7\glcondEnvLineTailCm{?}}

		\newcommand{\envHLineCSNmePdPt}[8][*]{\begin{eqnarray}#3&#2&#4,\\\glcdH#5&#2&#6,\\\glcdH#7&#2&#8\glcondEnvLineTailPd{?}}
		\newcommand{\envHLineCSCmNmePt}[8][*]{\begin{eqnarray}#3&#2&#4,\\\glcdH#5&#2&#6,\\\glcdH#7&#2&#8\glcondEnvLineTailCm{?}}
		\newcommand{\envHLineCSNmePdPte}[9][*]{\begin{eqnarray}#2&#3&#4,\\\glcdH#5&#6&#7,\\\glcdH#8&#2&#9\glcondEnvLineTailPd{?}}
		\newcommand{\envHLineCSCmNmePte}[9][*]{\begin{eqnarray}#2&#3&#4,\\\glcdH#5&#6&#7,\\\glcdH#8&#2&#9\glcondEnvLineTailCm{?}}
		
		\newcommand{\envHLineCECm}[9][*]{\glcondEnvLineHead{#1} #2&=&#3,\nonumber\\\glcdH#4&=&#5,\nonumber\\\glcdH#6&=&#7,\\\glcdH#8&=&#9\nonumber\glcondEnvLineTailCm{#1}}
		
		\newcommand{\envHLineCENmePd}[9][*]{\begin{eqnarray} #2&=&#3,\\\glcdH#4&=&#5,\\\glcdH#6&=&#7,\\\glcdH#8&=&#9\glcondEnvLineTailPd{?}}
		\newcommand{\envHLineCEDefNmePd}[9][*]{\begin{eqnarray} #2&:=&#3,\\\glcdH#4&:=&#5,\\\glcdH#6&:=&#7,\\\glcdH#8&:=&#9\glcondEnvLineTailPd{?}}
		\newcommand{\envHLineCECmNme}[9][*]{\begin{eqnarray} #2&=&#3,\\\glcdH#4&=&#5,\\\glcdH#6&=&#7,\\\glcdH#8&=&#9\glcondEnvLineTailCm{?}}
		\newcommand{\envHLineCECmDefNme}[9][*]{\begin{eqnarray} #2&:=&#3,\\\glcdH#4&:=&#5,\\\glcdH#6&:=&#7,\\\glcdH#8&:=&#9\glcondEnvLineTailCm{?}}

			\newcommand{\pccondPaOMathEnvCmPdPar}[1]{\ifx#1p \\\glcdH \fi}

		\newcommand{\envPLinePd}[2][*]{\glcondEnvLineHead{#1} #2\glcondEnvLineTailPd{#1}}
		\newcommand{\envPLineCm}[2][*]{\glcondEnvLineHead{#1} #2\glcondEnvLineTailCm{#1}}
		\newcommand{\envOTLineCm}[4][*]{\glcondEnvLineHead{#1} #2\lnAP{=}#3\lnP{=}#4,\glcondEnvLineTail{#1}}
		\newcommand{\envOTLinePd}[4][*]{\glcondEnvLineHead{#1} #2\lnAP{=}#3\lnP{=}#4.\glcondEnvLineTail{#1}}
		\newcommand{\envOTLineCmDef}[4][*]{\glcondEnvLineHead{#1} #2\lnAP{:=}#3\lnP{=}#4,\glcondEnvLineTail{#1}}
		
		\newcommand{\envOTLineCmPt}[5][*]{\glcondEnvLineHead{#1} #3\lnAP{#2}#4\lnP{#2}#5\glcondEnvLineTailCm{#1}}
		\newcommand{\envOTLinePdPt}[5][*]{\glcondEnvLineHead{#1} #3\lnAP{#2}#4\lnP{#2}#5\glcondEnvLineTailPd{#1}}
			\newcommand{\envOTLineThCm}{\envOTLineCm}

		\newcommand{\envOFLineCm}[5][*]{\glcondEnvLineHead{#1} #2\lnAP{=}#3\lnP{=}#4\lnP{=}#5,\glcondEnvLineTail{#1}}
		\newcommand{\envOFLinePd}[5][*]{\glcondEnvLineHead{#1} #2\lnAP{=}#3\lnP{=}#4\lnP{=}#5.\glcondEnvLineTail{#1}}

		\newcommand{\envMOCm}[2][*]{$\ifx d#1 \displaystyle \fi#2$,}
		\newcommand{\envMOPd}[2][*]{$\ifx d#1 \displaystyle \fi#2$.}
		
		\newcommand{\envMTCm}[3][*]{$\ifx d#1 \displaystyle \fi#2=#3$,}
		\newcommand{\envMTPd}[3][*]{$\ifx d#1 \displaystyle \fi#2=#3$.}
		\newcommand{\envMTCmDef}[3][*]{$\ifx d#1 \displaystyle \fi#2:=#3$,}
		\newcommand{\envMTDefPd}[3][*]{$\ifx d#1 \displaystyle \fi#2:=#3$.}
		\newcommand{\envMTCmPt}[4][*]{$\ifx d#1 \displaystyle \fi#3#2#4$,}
		\newcommand{\envMTPdPt}[4][*]{$\ifx d#1 \displaystyle \fi#3#2#4$.}
			\newcommand{\envMCm}{\envMTCm}
			\newcommand{\envMPd}{\envMTPd}

			\newcommand{\envMCmPt}{\envMTCmPt}
			
		\newcommand{\envMThCm}[4][*]{$\ifx d#1 \displaystyle \fi#2=#3=#4$,}
		\newcommand{\envMThPd}[4][*]{$\ifx d#1 \displaystyle \fi#2=#3=#4$.}
		\newcommand{\envMThCmPt}[5][*]{$\ifx d#1 \displaystyle \fi#3#2#4#2#5$,}
		\newcommand{\envMThPdPt}[5][*]{$\ifx d#1 \displaystyle \fi#3#2#4#2#5$.}
		\newcommand{\envMFCm}[5][*]{$\ifx d#1 \displaystyle \fi#2=#3=#4=#5$,}
		\newcommand{\envMFPd}[5][*]{$\ifx d#1 \displaystyle \fi#2=#3=#4=#5$.}
		\newcommand{\envMFCmPt}[6][*]{$\ifx d#1 \displaystyle \fi#3#2#4#2#5#2#6$,}
		\newcommand{\envMFPdPt}[6][*]{$\ifx d#1 \displaystyle \fi#3#2#4#2#5#2#6$.}

		\newcommand{\envHLineCFCmNm}[5][*]{ \begin{equation}\begin{split} \ifx*#1 \text{[ERROR;need label name]} \else \label{#1} \fi #2&\lnP{=}#3,\\#4&\lnP{=}#5, \end{split}\end{equation} }
		\newcommand{\envHLineCFNm}[5][*]{ \begin{equation}\begin{split} \ifx*#1 \text{[ERROR;need label name]} \else \label{#1} \fi #2&\lnP{=}#3\\#4&\lnP{=}#5, \end{split}\end{equation} }
		\newcommand{\envHLineCFNmPd}[5][*]{ \begin{equation}\begin{split} \ifx*#1 \text{[ERROR;need label name]} \else \label{#1} \fi #2&\lnP{=}#3,\\#4&\lnP{=}#5. \end{split}\end{equation} }
		\newcommand{\envHLineCFCmDefNm}[5][*]{ \begin{equation}\begin{split} \ifx*#1 \text{[ERROR;need label name]} \else \label{#1} \fi #2&\lnP{:=}#3,\\#4&\lnP{:=}#5, \end{split}\end{equation} }
		\newcommand{\envHLineCFDefNm}[5][*]{ \begin{equation}\begin{split} \ifx*#1 \text{[ERROR;need label name]} \else \label{#1} \fi #2&\lnP{:=}#3\\#4&\lnP{:=}#5, \end{split}\end{equation} }
		\newcommand{\envHLineCFDefNmPd}[5][*]{ \begin{equation}\begin{split} \ifx*#1 \text{[ERROR;need label name]} \else \label{#1} \fi #2&\lnP{:=}#3,\\#4&\lnP{:=}#5. \end{split}\end{equation} }
		\newcommand{\envHLineCSCmNm}[7][*]{ \begin{equation}\begin{split} \ifx*#1 \text{[ERROR;need label name]} \else \label{#1} \fi #2&\lnP{=}#3,\\#4&\lnP{=}#5,\\#6&\lnP{=}#7 \end{split}\end{equation} }
		\newcommand{\envHLineCSNm}[7][*]{ \begin{equation}\begin{split} \ifx*#1 \text{[ERROR;need label name]} \else \label{#1} \fi #2&\lnP{=}#3\\#4&\lnP{=}#5\\#6&\lnP{=}#7 \end{split}\end{equation} }
		\newcommand{\envHLineCSNmPd}[7][*]{ \begin{equation}\begin{split} \ifx*#1 \text{[ERROR;need label name]} \else \label{#1} \fi #2&\lnP{=}#3,\\#4&\lnP{=}#5,\\#6&\lnP{=}#7. \end{split}\end{equation} }
		\newcommand{\envHLineCSCmDefNm}[7][*]{ \begin{equation}\begin{split} \ifx*#1 \text{[ERROR;need label name]} \else \label{#1} \fi #2&\lnP{:=}#3,\\#4&\lnP{:=}#5,\\#6&\lnP{:=}#7 \end{split}\end{equation} }
		\newcommand{\envHLineCSDefNm}[7][*]{ \begin{equation}\begin{split} \ifx*#1 \text{[ERROR;need label name]} \else \label{#1} \fi #2&\lnP{:=}#3\\#4&\lnP{:=}#5\\#6&\lnP{:=}#7 \end{split}\end{equation} }
		\newcommand{\envHLineCSDefNmPd}[7][*]{ \begin{equation}\begin{split} \ifx*#1 \text{[ERROR;need label name]} \else \label{#1} \fi #2&\lnP{:=}#3,\\#4&\lnP{:=}#5,\\#6&\lnP{:=}#7. \end{split}\end{equation} }
		\newcommand{\envHLineCECmNm}[9][*]{ \begin{equation}\begin{split} \ifx*#1 \text{[ERROR;need label name]} \else \label{#1} \fi #2&\lnP{=}#3,\\#4&\lnP{=}#5,\\#6&\lnP{=}#7,\\#8&\lnP{=}#9,  \end{split}\end{equation} }
		\newcommand{\envHLineCENm}[9][*]{ \begin{equation}\begin{split} \ifx*#1 \text{[ERROR;need label name]} \else \label{#1} \fi #2&\lnP{=}#3\\#4&\lnP{=}#5\\#6&\lnP{=}#7\\#8&\lnP{=}#9  \end{split}\end{equation} }
		\newcommand{\envHLineCENmPd}[9][*]{ \begin{equation}\begin{split} \ifx*#1 \text{[ERROR;need label name]} \else \label{#1} \fi #2&\lnP{=}#3,\\#4&\lnP{=}#5,\\#6&\lnP{=}#7,\\#8&\lnP{=}#9.  \end{split}\end{equation} }
		\newcommand{\envHLineCECmDefNm}[9][*]{ \begin{equation}\begin{split} \ifx*#1 \text{[ERROR;need label name]} \else \label{#1} \fi #2&\lnP{:=}#3,\\#4&\lnP{:=}#5,\\#6&\lnP{:=}#7,\\#8&\lnP{:=}#9,  \end{split}\end{equation} }
		\newcommand{\envHLineCEDefNm}[9][*]{ \begin{equation}\begin{split} \ifx*#1 \text{[ERROR;need label name]} \else \label{#1} \fi #2&\lnP{:=}#3\\#4&\lnP{:=}#5\\#6&\lnP{:=}#7\\#8&\lnP{:=}#9  \end{split}\end{equation} }
		\newcommand{\envHLineCEDefNmPd}[9][*]{ \begin{equation}\begin{split} \ifx*#1 \text{[ERROR;need label name]} \else \label{#1} \fi #2&\lnP{:=}#3,\\#4&\lnP{:=}#5,\\#6&\lnP{:=}#7,\\#8&\lnP{:=}#9.  \end{split}\end{equation} }
	\newcommand{\envMLineTPd}[3][*]{ \ifx*#1 \begin{multline*} #2\lnP{=}#3.\end{multline*} \else \begin{multline} \label{#1} #2\lnP{=}#3.\end{multline} \fi }
	\newcommand{\envMLineTCm}[3][*]{ \ifx*#1 \begin{multline*} #2\lnP{=}#3,\end{multline*} \else \begin{multline} \label{#1} #2\lnP{=}#3,\end{multline} \fi }
	\newcommand{\envMLineTDefPd}[3][*]{ \ifx*#1 \begin{multline*} #2\lnP{:=}#3.\end{multline*} \else \begin{multline} \label{#1} #2\lnP{:=}#3.\end{multline} \fi }
	\newcommand{\envMLineTCmDef}[3][*]{ \ifx*#1 \begin{multline*} #2\lnP{:=}#3,\end{multline*} \else \begin{multline} \label{#1} #2\lnP{:=}#3,\end{multline} \fi }

	\newcommand{\envCaseTPd}[3][?]{\begin{cases} \glcondDis{#1}#2,\lcparaCase\\\glcondDis{#1}#3.\end{cases}}

	\newcommand{\envItemEmPa}[2][*]{\vspace{-5pt}\ifx*#1\begin{enumerate}\else\begin{enumerate}[#1]\fi \renewcommand{\itemsep}{0pt}\renewcommand{\itemindent}{0pt}\renewcommand{\theenumi}{\bfTx{(\roman{enumi})}}#2\end{enumerate}}

	\DeclareFontEncoding{OT2}{}{}
	\DeclareFontSubstitution{OT2}{cmr}{m}{n}
	\DeclareFontFamily{OT2}{cmr}{\hyphenchar\font45}
	\DeclareFontShape{OT2}{cmr}{m}{n}{<5><6><7><8><9>gen*wncyr <10><10.95><12><14.4><17.28><20.74><24.88>wncyr10}{}
	\DeclareFontShape{OT2}{cmr}{b}{n}{<5><6><7><8><9>gen*wncyb<10><10.95><12><14.4><17.28><20.74><24.88>wncyb10}{}
	\DeclareMathAlphabet{\mathcyr}{OT2}{cmr}{m}{n}
	\DeclareMathAlphabet{\mathcyb}{OT2}{cmr}{b}{n}
	\SetMathAlphabet{\mathcyr}{bold}{OT2}{cmr}{b}{n}

	\newcommand{\sh}{\mathcyr{sh}}
	\newcommand{\shS}[1][\;]{#1\mathcyr{sh}#1}

					\newcommand{\wspZV}[2][]{\mathcal{Z}_{#2}}			\newcommand{\wdspZV}[3][]{\mathcal{Z}_{#2}^{(#3)}}
\newcommand{\wspZVd}{\wdspZV}
					\newcommand{\wspPZV}[2][]{\mathcal{P}_{#2}}			

	\newcommand{\alHfnd}{\mathfrak{H}}		
	\newcommand{\alHR}{\alHfnd^1}		\newcommand{\alHRh}{\alHfnd_*^1}			\newcommand{\alHRs}{\alHfnd_\sh^1}		
	\newcommand{\alHN}{\alHfnd^0}		\newcommand{\alHNh}{\alHfnd_*^0}			\newcommand{\alHNs}{\alHfnd_\sh^0}	
	\newcommand{\alHF}{\alHfnd}							
		
	\newcommand{\fcZBlettHelp}{\bullet}
	\newcommand{\fcZBN}[1][?]{\zeta^\fcZBlettHelp}			
	\newcommand{\fcZHN}[1][?]{\zeta^*}				\newcommand{\fcZHO}[2][n]{\pwFc[#1]{\zeta}{*}{#2}}					\newcommand{\fcZH}{\fcZHO}
	\newcommand{\fcZRN}[1][?]{\zeta^\ddagger}		\newcommand{\fcZRO}[2][n]{\pwFc[#1]{\zeta}{\ddagger}{#2}}			\newcommand{\fcZR}{\fcZRO}			
	\newcommand{\fcZSN}[1][?]{\zeta^\sh}			\newcommand{\fcZSO}[2][n]{\pwFc[#1]{\zeta}{\sh}{#2}}				\newcommand{\fcZS}{\fcZSO}
								
	\newcommand{\LetHelpFcCh}{\chi}
	\newcommand{\fcChBN}[1][?]{\LetHelpFcCh^\fcZBlettHelp}			
	\newcommand{\fcChHN}[1][?]{\LetHelpFcCh^*}									
	\newcommand{\fcChRN}[1][?]{\LetHelpFcCh^\ddagger}								
	\newcommand{\fcChSN}[1][?]{\LetHelpFcCh^\sh}							

	\newcommand{\fcChBmN}[1][?]{\mLt{\LetHelpFcCh}^\fcZBlettHelp}			
	\newcommand{\fcChHmN}[1][?]{\mLt{\LetHelpFcCh}^*}							
	\newcommand{\fcChRmN}[1][?]{\mLt{\LetHelpFcCh}^\ddagger}								
	\newcommand{\fcChSmN}[1][?]{\mLt{\LetHelpFcCh}^\sh}							
								
	\newcommand{\mpEVlett}{Z}
	\newcommand{\mpEVN}[2][n]{\Fc[#1]{\mpEVlett}{#2}}				
	\newcommand{\mpEVH}[2][n]{\pwFc[#1]{\mpEVlett}{*}{#2}}		\newcommand{\mpEVHi}[3][n]{\ipFc[#1]{\mpEVlett}{#2}{*}{#3}}			
	\newcommand{\mpEVS}[2][n]{\pwFc[#1]{\mpEVlett}{\sh}{#2}}		\newcommand{\mpEVSi}[3][n]{\ipFc[#1]{\mpEVlett}{#2}{\sh}{#3}}

	\newcommand{\gfcZlett}{\mathfrak{Z}}		\newcommand{\ggfcZlett}{\mathfrak{\bLt{Z}}}	
				
		\newcommand{\gfcZH}[4][n]{\ipFc[#1]{\gfcZlett}{#2}{*(#3)}{#4}}								\newcommand{\gfcZHbb}[2][?]{\gfcZlett_#2^*}
					\newcommand{\gfcZRb}[3][n]{\ipFc[#1]{\gfcZlett}{#2}{\ddagger}{#3}}		\newcommand{\gfcZRbb}[2][?]{\gfcZlett_#2^\ddagger}
		\newcommand{\gfcZS}[4][n]{\ipFc[#1]{\gfcZlett}{#2}{\sh(#3)}{#4}}							\newcommand{\gfcZSbb}[2][?]{\gfcZlett_#2^\sh}
									
												\newcommand{\ggfcZbb}[1][?]{\ggfcZlett}
		\newcommand{\ggfcZH}[3][n]{\pwFc[#1]{\ggfcZlett}{*(#2)}{#3}}								\newcommand{\ggfcZHbb}[1][?]{\ggfcZlett^*}
							\newcommand{\ggfcZRbb}[1][?]{\ggfcZlett^\ddagger}
		\newcommand{\ggfcZS}[3][n]{\pwFc[#1]{\ggfcZlett}{\sh(#2)}{#3}}							\newcommand{\ggfcZSbb}[1][?]{\ggfcZlett^\sh}

	\newcommand{\lbfL}{{\bf l}}
	
\newcommand{\lgpeD}[2][?]{\ifx#1?\mathsf{D}_{#2}\else\mathsf{D}_{#2}^{(#1)}\fi}
\newcommand{\lgpeR}[1][?]{\ifx#1?\mathsf{R}\else\mathsf{R}^{(#1)}\fi}
\newcommand{\lgpeRO}[2][?]{\ifx#1?\mathsf{R}_{#2}\else\mathsf{R}_{#2}^{(#1)}\fi}
\newcommand{\lgpeRm}[1][?]{\ifx#1?\mLt{\mathsf{R}}\else\mLt{\mathsf{R}}^{(#1)}\fi}
\newcommand{\lgpeP}[1][?]{\ifx#1?\mathsf{P}\else\mathsf{P}^{(#1)}\fi}
\newcommand{\lgpePO}[2][?]{\ifx#1?\mathsf{P}_{#2}\else\mathsf{P}_{#2}^{(#1)}\fi}
\newcommand{\lgpePi}[1][?]{\ifx#1?\mathsf{P}^{-1}\else(\mathsf{P}^{(#1)})^{-1}\fi}
\newcommand{\lgpePOi}[2][?]{\ifx#1?(\mathsf{P}_{#2})^{-1}\else(\mathsf{P}_{#2}^{(#1)})^{-1}\fi}

	\newcommand{\lsstSh}[1][?]{\ifx#1?{\mathsf{S}}_{\sh}\else{\mathsf{S}}_{\sh}^{(#1)}\fi}
	\newcommand{\lsstShT}[2][?]{\ifx#1?{\mathsf{S}}_{\sh,#2}\else{\mathsf{S}}_{\sh,#2}^{(#1)}\fi}

	\renewcommand{\matu}[1][?]{\ifx#1?\mathsf{I}\else {\mathsf{I}^{(#1)}}\fi}
	\renewcommand{\gpu}[2][?]{\ifx ?#1e^{(#2)}\else e\fi}
	\newcommand{\popR}[2][?]{\bkR[#1]{#2}^\vee}
	\newcommand{\pbfL}[1][?]{\bfTx{l}}
	\newcommand{\pbfK}[1][?]{\bfTx{k}}
	\newcommand{\pbfX}[1][?]{\bfTx{x}}
	\newcommand{\pivIN}{i}		\newcommand{\pivIO}[2][n]{\Fc[#1]{\pivIN}{#2}}
	\newcommand{\pivTRN}{t}	\newcommand{\pivTRO}[2][n]{\Fc[#1]{\pivTRN}{#2}}
			
	\newcommand{\pmpRH}[2][n]{\Fc[#1]{\rho}{#2}}
	\newcommand{\pmpRG}[2][n]{\ipFc[#1]{\mathrm{reg}}{\sh}{T}{#2}}		\newcommand{\mpRG}{\pmpRG}
	\newcommand{\pmpRGb}[2][n]{\idFc[#1]{\mathrm{reg}}{\sh}{#2}}		
	\newcommand{\popSP}[2][?]{\bkR[#1]{#2}^\sharp}				\newcommand{\popSPr}[2][n]{(#2)^\sharp}
				\newcommand{\popSPrT}[3][n]{\pwFc[#1]{(#2)}{\sharp}{#3}}
	\newcommand{\pnmSN}{s}				\newcommand{\pnmS}[2][?]{\pnmSN_{#2}}
	\newcommand{\pnmSrN}{s^\ddagger}		\newcommand{\pnmSr}[2][?]{\pnmSrN_{#2}}
	\newcommand{\pnmRN}{r}				\newcommand{\pnmR}[2][?]{\pnmRN_{#2}}
	
	\newcommand{\pggfcZe}[3][n]{\ipFc[#1]{F}{#2,\ep}{\sh}{#3}}
	\renewcommand{\gpSym}[2][?]{S_{#2}}
	\newcommand{\lfgCYCO}[2][?]{C_{#2}}		\newcommand{\lfgCYCT}[3][?]{C_{#2}^{(#3)}}
	\newcommand{\lfgEPN}{\ep}		\newcommand{\lfgEPO}[2][?]{\lfgEPN^{(#2)}}				\newcommand{\lfgEP}{\lfgEPN}
	\newcommand{\lfgPN}{P}			\newcommand{\lfgPO}[2][?]{\lfgPN^{(#2)}}				\newcommand{\lfgP}{\lfgPN}
								\newcommand{\lfgPiO}[2][b]{\bkR[#1]{\lfgPN^{(#2)}}^{-1}}	\newcommand{\lfgPi}{\lfgPN^{-1}}
	\newcommand{\lfgSHiN}{sh}		\newcommand{\lfgSHiT}[3][?]{\lfgSHiN_{#2}^{(#3)}}	
								\newcommand{\lfgSHiO}[2][?]{\lfgSHiN_{#2}}		\newcommand{\lfgSHi}{\lfgSHiO}
	\newcommand{\lfgTN}{T}			\newcommand{\lfgTT}[3][?]{\lfgTN_{#2}^{(#3)}}	
								\newcommand{\lfgTO}[2][?]{\lfgTN_{#2}}					\newcommand{\lfgT}{\lfgTN}
	\newcommand{\lfgTUN}{\tau}		\newcommand{\lfgTUO}[2][?]{\tau_{#2}}	\newcommand{\lfgTUOO}[2][?]{\tau^{(#2)}}	
	\newcommand{\lfgUTN}{e}		\newcommand{\lfgUTO}[2][?]{\lfgUTN^{(#2)}}		\newcommand{\lfgUT}{\lfgUTN}
	\renewcommand{\gfcZlett}{f}		\renewcommand{\ggfcZlett}{F}	
											
		\renewcommand{\gfcZH}[4][n]{\ipFc[#1]{\gfcZlett}{#2,#3}{*}{#4}}								\renewcommand{\gfcZHbb}[2][?]{\gfcZlett_#2^*}
					\renewcommand{\gfcZRb}[3][n]{\ipFc[#1]{\gfcZlett}{#2}{\ddagger}{#3}}	\renewcommand{\gfcZRbb}[2][?]{\gfcZlett_#2^\ddagger}
		\renewcommand{\gfcZS}[4][n]{\ipFc[#1]{\gfcZlett}{#2,#3}{\sh}{#4}}							\renewcommand{\gfcZSbb}[2][?]{\gfcZlett_#2^\sh}
											\renewcommand{\ggfcZbb}[1][?]{\ggfcZlett}
		\renewcommand{\ggfcZH}[3][n]{\ipFc[#1]{\ggfcZlett}{#2}{*}{#3}}								\renewcommand{\ggfcZHbb}[1][?]{\ggfcZlett^*}
						\renewcommand{\ggfcZRbb}[1][?]{\ggfcZlett^\ddagger}
		\renewcommand{\ggfcZS}[3][n]{\ipFc[#1]{\ggfcZlett}{#2}{\sh}{#3}}							\renewcommand{\ggfcZSbb}[1][?]{\ggfcZlett^\sh}

\allowdisplaybreaks[4]
	\title{\mainTitle}
	\author{\authorName
			\thanks{\organizationNameFst, \placeAddressFst}
			\mbox{}
			\thanks{\organizationNameScd, \departmentNameScd, \placeAddressScd}
		}
	\date{}

\begin{document}
\maketitle
\renewcommand{\thefootnote}{\fnsymbol{footnote}}
\footnote[0]{e-mail : \emailAddressFst}
\renewcommand{\thefootnote}{\arabic{footnote}}\setcounter{footnote}{0}
\vPack[30]

\begin{abstract}
Riemann zeta values are generalized to multiple zeta values (MZVs) by use of nested sums,
	and MZVs are generalized to regularized multiple zeta values (RMZVs) by regularization of divergent infinite series. 
In the present paper,
	we prove congruence identities of RMZVs of depth $n$ involving a pair of index sets;
	the congruence relation is given by the vector space spanned by MZVs of depth $n-1$ and products of MZVs.
We also obtain a proof of the parity result,
	and a congruence sum formula for MZVs.
\end{abstract}

\section{Introduction and statement of results} \label{sectOne}
The special values of the Riemann zeta function
	\envM{ \fcZ{s} }{ \SmT{m=1}{\infty} \opF[s]{1}{m^s} } at integer arguments have attracted many mathematicians including Euler, 
	who solved the Basel problem,
	$\fcZ{2}=\opF[s]{\pi^2}{6}$.
These values play an important role in number theory,
	and have applications in many areas \cite{Zagier94}.
A multiple zeta value (MZV)
	\envHLine
	{
		\fcZ{l_1,\ldots,l_n}
	}
	{
		\Sm{m_1>\cdots>m_n>0} \opF{1}{ \pw{m_1}{l_1}\cdots\pw{m_n}{l_n} }
	} 
	is a generalization of such a special value,
	and defined for an (ordered) index set  $(l_1,\ldots,l_n)$ of positive integers with $l_1\geq2$.
The condition $l_1\geq2$ ensures the convergence.
We call $l=l_1+\cdots+l_n$ and $n$ the weight and depth, respectively.
MZVs were first studied by Euler \cite{Euler1775} in the case of depth $2$.

Let $\fcZH{l_1,\ldots,l_n}$ and $\fcZS{l_1,\ldots,l_n}$ be regularized multiple zeta values (RMZVs) of harmonic and shuffle types, 
	which are defined in \cite{IKZ06} as generalizations of MZVs. 
(The details of these values will be introduced in \refSec{sectTwoOne}.)
Both types of RMZVs are MZVs if $l_1\geq2$,
	but they are also defined for index sets $(l_1,\ldots,l_n)$ with $l_1=1$, unlike MZVs.
We know from \cite[(2,3)]{IKZ06} that
	the difference of $\fcZH{l_1,\ldots,l_n}$ and $\fcZS{l_1,\ldots,l_n}$
	can be expressed as a linear combination in the vector space $\wspPZV{l}$ over $\setQ$,
	which is spanned by products of MZVs with total weight $l$,
	where the products include the single zeta value $\fcZ{l}$.
That is,
	$\fcZH{l_1,\ldots,l_n} \equiv \fcZS{l_1,\ldots,l_n} \ \modd \wspPZV{l}$ (see \refEq{2.1_Prop1_Eq1} in \refSec{sectTwoOne}),
	and so we may use the same symbol $\fcZR[]{}$ for both $\fcZH[]{}$ and $\fcZS[]{}$
	when we consider a congruence relation modulo a vector space including $\wspPZV{l}$.

MZVs satisfy an interesting property called the parity result,
	which was proved in \cite{IKZ06,Tsumura04} for general depth $n$ independently (see \cite{BBG95,Euler1775,Tornheim50} for depth $2$ and \cite{BG96} for depth $3$),
	and plays an important role in the study of the depth-graded algebra of MZVs \cite{Brown13Ax,IKZ06}.
The parity result states a reducibility of MZVs when the weight and depth have opposite parity.
Since every RMZV is written in terms of MZVs of depth $n$ and products of MZVs with total weight $l$ (see \refEq{2.1_Prop1_Eq2} in \refSec{sectTwoOne}),
	RMZVs also satisfy the parity result.
Thus,
	\envHLineCmPt[1_PL_EqqParityResult]{\equiv}
	{
		\fcZR{l_1,\ldots,l_n}
	}
	{
		0
		\qquad
		\modd \wspZVd{l}{<n} + \wspPZV{l}
		\qquad
		(\text{$l+n$ is odd})
	}
	where $\wspZVd{l}{<n}$ is the vector space over $\setQ$ spanned by MZVs of weight $l$ and depth less than $n$.

For an index set $\lbfL=(l_1,l_2,\ldots,l_n)$,
	we define its reverse $\popR{\lbfL}=\popR{(l_1,l_2,\ldots,l_n)}$ by
	\envHLineDefPd
	{
		\popR{(l_1,l_2,\ldots,l_n)}
	}
	{
		(l_n,l_{n-1},\ldots,l_1)
	}
In the proof of \refEq{1_PL_EqqParityResult},
	Ihara, Kaneko, and Zagier \cite{IKZ06} showed the following congruence identities of RMZVs involving the pair $(\lbfL,\popR{\lbfL})$ of index sets:
	\begin{align}
		\fcZR{\lbfL}	&\equiv\mo^{l-1} \fcZR{\popR{\lbfL}}				&\hspace{-80pt}\modd \wspZVd{l}{<n} + \wspPZV{l}, 		\label{1_PL_EqqRefRsltFromIKZ1}\\
		\fcZR{\lbfL}	&\equiv\mo^{n-1} \fcZR{\popR{\lbfL}}				&\hspace{-80pt}\modd \wspZVd{l}{<n} + \wspPZV{l}. 		\label{1_PL_EqqRefRsltFromIKZ2}
	\end{align}
Identity \refEq{1_PL_EqqRefRsltFromIKZ1} does not explicitly appear in \cite{IKZ06}, 
	but it is easily proved
	because \refEq{1_PL_EqqParityResult} and \refEq{1_PL_EqqRefRsltFromIKZ2} give \refEq{1_PL_EqqRefRsltFromIKZ1} and \refEq{1_PL_EqqRefRsltFromIKZ2},
	and vice versa.
Identity \refEq{1_PL_EqqRefRsltFromIKZ2} is given in \cite[(8.6)]{IKZ06}.

In the present paper,
	we improve \refEq{1_PL_EqqRefRsltFromIKZ1} 
	by reducing the vector space giving the congruence relation $\equiv$ from $\wspZVd{l}{<n} + \wspPZV{l}$ to its subspace,
	and prove two congruence identities of RMZVs involving a pair of index sets modulo the subspace.
The subspace is given by $\wspZVd{l}{n-1} + \wspPZV{l}$,
	where $\wspZVd{l}{d}$ is the vector space spanned by MZVs of weight $l$ and depth $d$.
Note that $\wspZVd{l}{n-1} \subset\wspZVd{l}{<n}=\SmT{d=1}{n-1}\wspZVd{l}{d}$.

Our results are stated as follows.
\begin{theorem}\label{1_Thm1}
Identity \refEq{1_PL_EqqRefRsltFromIKZ1} holds 
	even if the vector space $\wspZVd{l}{<n} + \wspPZV{l}$ is changed to the subspace $\wspZVd{l}{n-1} + \wspPZV{l}$.
In particular, 
	we obtain another proof of the parity result \refEq{1_PL_EqqParityResult}.
\end{theorem}
\begin{theorem}\label{1_Thm2}
We define two index sets of positive integers by $\pbfK_i=(l_1,\ldots,l_{i-1})$ and $\pbfL_i=(l_{i+1},\ldots,l_n)$ for an integer $i$ with $1\leq i \leq n$.
Then we have
	\begin{align}
		\fcZR{\pbfK_i,\osTx{1}{i\,th},\pbfL_i}	&\equiv\fcZR{\pbfL_i,\osTx{1}{j\,th},\pbfK_i}						&\hspace{-50pt}\modd \wspZVd{l}{n-1} + \wspPZV{l}, 	
			\label{1_Thm2_EqqRefRslt1}\\
		\fcZR{\pbfK_i,\osTx{1}{i\,th},\pbfL_i}	&\equiv\mo^{l-1} \fcZR{ \popR{\pbfK_i},\osTx{1}{i\,th},\popR{\pbfL_i}}		&\hspace{-50pt}\modd \wspZVd{l}{n-1} + \wspPZV{l},
			\label{1_Thm2_EqqRefRslt2}
	\end{align}	
	where $\pbfK_1$ and $\pbfL_n$ mean the empty set, and $j=n+1-i$.
\end{theorem}

As a corollary of \refThm{1_Thm2},
	we can obtain the following congruence sum formula for MZVs with the help of the restricted sum formula \cite{ELO09} and Ohno's relations \cite{Ohno99}.
\begin{corollary}\label{1_Cor1}
Let $l,m,n$ be integers with $l\geq m+n$ and $m, n\geq1$. 
If weight $l$ is even,
	\envHLineCmPt[1_Cor1_EqqRestrictedOhnoRel]{\equiv}
	{
		\tpSm{l_1\geq m+1, l_2,\ldots,l_n\geq1}{l_1+\cdots+l_n=l} \fcZ{l_1,\ldots,l_n}
	}
	{
		0
		\qquad
		\modd\wspZVd{l}{d+1} + \wspPZV{l}
	}
	where $d\in\SetO{m,n}$.
In particular,
	the left-hand side of \refEq{1_Cor1_EqqRestrictedOhnoRel} is congruent to $0$ modulo $\wspZVd{l}{<n} + \wspPZV{l}$ under the additional condition $m+1< n$.
\end{corollary}

The present paper is organized as follows.
We devote \refSect{sectTwo} to preliminaries;
	we review the algebraic setup of RMZVs and some of their properties by referring to \cite{Hoffman97,IKZ06} in \refSect{sectTwoOne},
	and we show key identities in $\setZ[\gpGL{n}{\setZ}]$  by using results of \cite{IKZ06} in \refSect{sectTwoTwo},
	where $\gpGL{n}{\setZ}$ is the general linear group of degree $n$ over $\setZ$ and $\setZ[\gpGL{n}{\setZ}]$ is its group ring.
In \refSect{sectThree}, 
	we prove \refThm[s]{1_Thm1}, \ref{1_Thm2} and \refCor{1_Cor1},
	and give examples of \refThm{1_Thm2} for smaller weight.

\section{Preliminaries} \label{sectTwo}
\subsection{RMZVs of harmonic and shuffle types} \label{sectTwoOne}
Let $\alHF=\setQ\bkA{x,y}$ be the non-commutative polynomial algebra over $\setQ$ in two indeterminates $x$ and $y$,
	and $\alHN$ and $\alHR$ its subalgebras $\setQ+x\alHF y$ and $\setQ+\alHF y$, respectively.
These satisfy the inclusion relations $\alHN\subset\alHR\subset\alHF$.
Let $z_l$ denote $x^{l-1}y$ for any integer $l\geq1$.
Every word $w=w_0y$ in the set $\Set{x,y}$ with terminal letter $y$
	is expressed as $w=z_{l_1}\cdots z_{l_n}$ uniquely,
	and so $\alHR$ is the free algebra generated by $z_l\;(l=1,2,3,\ldots)$.
We define the harmonic product $*$ on $\alHR$ inductively by
	\envPLineTNme
	{\label{2_PL_DefHarProdUnit}&
		1 * w \lnP{=} w * 1 \lnP{=} w
	,&}
	{\label{2_PL_DefHarProdMain}&
		z_kw_1 * z_lw_2 \lnP{=} z_k(w_1*z_lw_2) + z_l(z_kw_1*w_2) + z_{k+l}(w_1*w_2)
	,&}
	for any integers $k,l\geq1$ and words $w, w_1, w_2\in\alHR$,
	and then extending by $\setQ$-bilinearity.
This product gives $\alHN$ and $\alHR$ structures of commutative $\setQ$-algebras \cite{Hoffman97},
	which we denote by $\alHNh$ and $\alHRh$, respectively;
	note that obviously $\alHNh$ is a subalgebra of $\alHRh$.
We also define the shuffle product $\sh$ on $\alHR$ inductively by
	\envPLineTNme
	{\label{2_PL_DefShuProdUnit}&
		1 \shS u \lnP{=} u \shS 1 \lnP{=} u
	,&}
	{\label{2_PL_DefShuProdMain}&
		uw_1\shS vw_2 \lnP{=} u(w_1 \shS vw_2) + v(uw_1 \shS w_2)
	,&}
	for any words $w, w_1, w_2\in\alHR$ and $u,v\in\Set{x,y}$, 
	and again extending by $\setQ$-bilinearity.
This product gives $\alHN$ and $\alHR$ the structures of commutative $\setQ$-algebras \cite{Reutenauer93},
	which we denote by $\alHNs$ and $\alHRs$, respectively;
	$\alHNs$ is a subalgebra of $\alHRs$.

Let $\mpEVN[]{}$ be the $\setQ$-linear map (called the evaluation map) from $\alHN$ to $\setR$ given by
	\envHLinePd
	{
		\mpEVN{z_{l_1}\cdots z_{l_n}}
	}
	{
		\fcZ{l_1,\ldots,l_n}
		\qquad
		(z_{l_1}\cdots z_{l_n}\in\alHN)
	}
We know from \cite{Hoffman97,Reutenauer93} that 
	$\mpEVN[]{}$ is homomorphic on both operations $*$ and $\sh$, 
	or satisfies
	\envOTLine
	{
		\mpEVN{w_1*w_2}
	}
	{
		\mpEVN{w_1 \shS w_2}
	}
	{
		\mpEVN{w_1} \mpEVN{w_2}
	}
	for $w_1,w_2\in\alHN$. 
Note that
	the weight and depth of $\fcZ{\lbfL_n}=\fcZ{l_1,\ldots,l_n}$ correspond to 
	the degree of $z_{l_1}\cdots z_{l_n}$ and the degree of $z_{l_1}\cdots z_{l_n}$ with respect to $y$,
	respectively,
	where we assign $x$ and $y$ degree $1$.
	
Let $\setR[T]$ be the polynomial ring in a single indeterminate with real coefficients.
Through the isomorphisms $\alHRh\simeq\alHNh[y]$ and $\alHRs\simeq\alHNs[y]$ proved in \cite{Hoffman97} and \cite{Reutenauer93},
	\etalTx{Ihara} \cite[\refPropA{1}]{IKZ06} considered 
	the algebra homomorphisms $\mpEVH[]{} : \alHRh \to \setR[T]$ and $\mpEVS[]{} : \alHRs \to \setR[T]$,
	respectively,
	which are uniquely characterized by the properties that it extends the evaluation map $\mpEVN[]{}$ and sends $y$ to $T$.		
For any word $w=z_{l_1}\cdots z_{l_n}\in\alHR$,
	we denote by $\mpEVHi{\lbfL_n}{T}$ and $\mpEVSi{\lbfL_n}{T}$ 
	the images under the maps $\mpEVH[]{}$ and $\mpEVS[]{}$ of the word $w$,
	and then RMZVs $\fcZH{\lbfL_n}$ and $\fcZS{\lbfL_n}$ of harmonic and shuffle types are defined by
	\envHLineCFLaaCm[2_PL_DefRMZV]
	{
		\fcZH{l_1,\ldots,l_n}
	}
	{
		\mpEVHi{l_1,\ldots,l_n}{0}
	}
	{
		\fcZS{l_1,\ldots,l_n}
	}
	{
		\mpEVSi{l_1,\ldots,l_n}{0}
	}
	respectively.
For example,
	we see from \refEq{2_PL_DefHarProdUnit}, \refEq{2_PL_DefHarProdMain}, \refEq{2_PL_DefShuProdUnit}, and \refEq{2_PL_DefShuProdMain}  that
	\envHLineCECm[2_PL_ExMapEV]
	{
		\mpEVHi{1}{T}
	}
	{
		\mpEVH{y}
	\lnP{=}
		T
	}
	{
		\mpEVSi{1}{T}
	}
	{
		\mpEVS{y}
	\lnP{=}
		T
	}
	{
		\mpEVHi{1,1}{T}
	}
	{
		\mpEVH[B]{\opF{1}{2}z_1*z_1-\opF{1}{2}z_2}
	\lnP{=}
		\opF{1}{2}T^2 - \opF{1}{2} \fcZ{2}
	}
	{
		\mpEVSi{1,1}{T}
	}
	{
		\mpEVS[B]{\opF{1}{2}y \shS y }
	\lnP{=}
		\opF{1}{2}T^2 
	}
	and so we have
	\envHLineCFLqqPd[2_PL_ExRMZV]
	{
		\fcZH{1}
	}
	{
		\fcZS{1}
	\lnP{=}
		\fcZS{1,1}
	\lnP{=}
		0
	}
	{
		\fcZH{1,1}
	}
	{
		-\opF{1}{2}\fcZ{2}
	}

We can derive the following properties for RMZVs from results proved in \cite{IKZ06}.

\begin{proposition}\label{2.1_Prop1}
Let $l_1,\ldots,l_n,l$ be positive integers with $l_1+\cdots+l_n=l$.
We have
	\begin{align}
		\fcZH{l_1,\ldots,l_n}	&\equiv\fcZS{l_1,\ldots,l_n}					&&\hspace{-30pt}	\modd \wspPZV{l}, 					\label{2.1_Prop1_Eq1}\\
		\fcZH{l_1,\ldots,l_n}	&\equiv\fcZS{l_1,\ldots,l_n}	\equiv	0		&&\hspace{-30pt}	\modd \wspZVd{l}{n} + \wspPZV{l}. 		\label{2.1_Prop1_Eq2}
	\end{align}
\end{proposition}
\envProof{
Let $\pmpRH[]{} : \setR[T] \to \setR[T]$ be the $\setR$-linear map which is determined by
	\envHLine[2.2_Prop1Proof_DefFuncRH]
	{
		\opF{ \pmpRH{ T^j } }{j!} 
	}
	{
		\SmT{i=0}{j} \gam_i \opF{T^{j-i}}{(j-i)!}
		\qquad
		(j=0,1,2,\ldots,)
	}
	and the $\setR$-linearity,
	where the coefficients $\gam_0=1$, $\gam_1=0$, $\gam_2=\opF[s]{\fcZ{2}}{2}$, 
	$\gam_3=-\opF[s]{\fcZ{3}}{3}$, $\gam_4=\opF[s]{(2\fcZ{4}+\fcZ{2}^2)}{8}$, \ldots\mbox{}
	are given by the generating function
	\envHLinePd
	{
		\exx[g]{\SmT{m=2}{\infty} \opF{\mo^m\fcZ{m}}{m} u^m }
	}
	{
		\SmT{i=0}{\infty} \gam_i u^i
	}
When $i\geq2$,
	we easily see that the coefficient $\gam_i$ can be written in terms of products of values $\fcZ{m}$ with total weight $i$,
	which yields the congruence equation 
	\envHLineCmPt[2.2_Prop1Proof_Eq1.1]{\equiv}
	{
		\fcZ{k_1,\ldots,k_d} \gam_i
	}
	{
		0
		\qquad
		\modd	\wspPZV{l}
	}
	where $\fcZ{k_1,\ldots,k_d}$ is an MZV of weight $l-i$.
\etalTx{Ihara} pointed out in \cite[\refSecA{2}]{IKZ06} that the coefficient $c_j^*$ of $T^j$ in $\mpEVHi{\lbfL_n}{T}$ is a $\setQ$-linear combination of MZVs of weight $l-j$, 
	that is,
	\envHLineCmPt[2.2_Prop1Proof_Eq1.2IKZ]{\in}
	{
		c_j^*
	}
	{
		\wspZV{l-j}
	}
	and also showed in \cite[\refThmA{1}]{IKZ06} that
	\envHLinePd[2.2_Prop1Proof_Eq1.3IKZ]
	{
		\pmpRH{\mpEVHi{l_1,\ldots,l_n}{T}}
	}
	{
		\mpEVSi{l_1,\ldots,l_n}{T}
	}
Therefore,
	because of \refEq{2.2_Prop1Proof_DefFuncRH}, \refEq{2.2_Prop1Proof_Eq1.1}, and \refEq{2.2_Prop1Proof_Eq1.2IKZ},
	substituting $T=0$ into \refEq{2.2_Prop1Proof_Eq1.3IKZ} gives \refEq{2.1_Prop1_Eq1}.
	
Let $\mpRG[]{}:\alHRs\to\alHNs[T]$ be the isomorphism which is uniquely characterized by the properties that it is the identity on $\alHNs$ and maps $y$ to $T$.
Note that
	$\mpEVS[]{} = \mpEVN[]{}\otimes \mpI[{\setQ[T]}] \circ \mpRG[]{}$ under the identification $ \alHNs \otimes \setQ[T] = \alHNs[T]$ (see \cite[\refSecA{3}]{IKZ06}).
Let $w=z_{l_1}\ldots z_{l_n}$ be a word in $\alHR$.
By the definitions of $\alHN$ and $\alHR$,
	there exist an integer $m\geq0$ and a word $w_0=xw_0'\in\alHN$ such that $w=y^mw_0$.
Then \etalTx{Ihara} \cite[\refPropA{8}]{IKZ06} showed that
	\envHLineCm
	{
		\pmpRGb{w}
	}
	{
		\mo^m x(y\shS w_0')
	}
	where $\pmpRGb[]{}$ is the homomorphism obtained by specializing $\mpRG[]{}$ to $T=0$.
By applying the map $\mpEVN[]{}$ to both sides of this equation,
	we obtain
	\envHLinePd[2.2_Prop1Proof_Eq2.1]
	{
		\fcZS{l_1,\ldots,l_n}
	}
	{
		\mo^m \mpEVN{ x(y^m\shS w_0') }
	}
For each word $v$ appearing in the expansion of $x(y^m\shS w_0')$,
	the degrees of $v$ and $xy^mw_0'$ 
	\respTx{the degrees of $v$ and $xy^mw_0'$ with respect to $y$} are same by the definition \refEq{2_PL_DefShuProdMain} of the operation $\sh$.
The degrees of $xy^mw_0'$ and $y^mxw_0'$ \respTx{the degrees of $xy^mw_0'$ and $y^mxw_0'$ with respect to $y$} are also same.
Since $w=y^mxw_0'$ and since the degree of $w$ and the degree of $w$ with respect to $y$ are $l=l_1+\cdots+l_n$ and $n$,
	respectively,
	the image of $x(y^m\shS w_0')$ under $\mpEVN[]{}$ is expressed in terms of MZVs of weight $l$ and depth $n$,
	i.e.,
	\envHLinePdPt[2.2_Prop1Proof_Eq2.2]{\equiv}
	{
		\mpEVN{ x(y^m\shS w_0') }
	}
	{
		0
		\qquad	\modd	\wspZVd{l}{n}
	}
Equating \refEq{2.2_Prop1Proof_Eq2.1} and \refEq{2.2_Prop1Proof_Eq2.2},
	we obtain
	\envHLineCm
	{
		\fcZS{l_1,\ldots,l_n}
	}
	{
		0
		\mod	\wspZVd{l}{n}
	}
	which together with \refEq{2.1_Prop1_Eq1} proves  \refEq{2.1_Prop1_Eq2}.
}

\subsection{Key identities in the group ring $\setZ[\gpGL{n}{\setZ}]$} \label{sectTwoTwo}
We begin by reviewing some notation about the group ring $\setZ[\gpGL{n}{\setZ}]$ introduced in \cite[\refSecA{8}]{IKZ06}.
Let $\gpSym{n}$ be the symmetric group of degree $n$
	and $\lfgUT=\lfgUTO{n}$ be its unit element.
We identify the elements of $\gpSym{n}$ with the permutation matrices in $\setZ[\gpGL{n}{\setZ}]$ in the usual way,
	that is,
	$\sig\in\gpSym{n}$ is identified with $(\delta_{i\sig(j)})_{1\leq i,j\leq n}\in\gpGL{n}{\setZ}$,
	where $\delta_{ij}$ is the Kronecker delta function.
We define the elements $\lfgEP=\lfgEPO{n}$ and $\lfgTO{j}=\lfgTT{j}{n} (j=0,\ldots,n-1)$ of order $2$ by
	\envHLineCFLaCm[2.2_PL_DefMinusAndReflect]
	{
		\lfgEPO{n}
	}
	{
		\MatTh{-1&&}{ &\raMTx{0}{\ddots}  &  }{&&-1}
	}
	{
		\lfgTT{j}{n}
	}
	{
		\MatT{1&\ldots&j&j+1&\ldots&n}{1&\ldots&j&n&\ldots&j+1}
	}
	respectively.
We set 
	\envM{
		\lfgTN
	}{
		\lfgTT{0}{n}
	}
	particularly, 
	which is identified with an anti-diagonal matrix,
	\envHLinePd
	{
		\lfgTN
	}
	{
		\MatTh{&&1}{ &\!\!\raMTx{-4}{\roMTx{85}{\ddots}}  &  }{1&&}
	}
We denote by $\lfgP=\lfgPO{n}$ the matrix given by
	\envHLineCFLqqPd[2.2_PL_DefP]
	{
		\lfgPO{n}
	}
	{
		\MatF{1&& &}{ -1&1&&}{ & \roMTx{-7}{\ddots}&\roMTx{-7}{\ddots}&  }{& & -1&1}
	}
	{
		\lfgPiO{n}
	}
	{
		\MatF{1&& &}{ 1&\,1&&}{ \vdots&\vdots&\roMTx{-7}{\ddots}&  }{1&\,1& \cdots&1}
	}
For each integer $j$ with $1\leq j \leq n-1$,
	we define the $j$-th shuffle element $\lfgSHi{j}=\lfgSHiT{j}{n}$ in $\setZ[\gpSym{n}]$ by
	\envHLineDefPd[2.2_PL_DefShuffle]
	{
		\lfgSHiT{j}{n}
	}
	{
		\tpSm{\sig\in\gpSym{n}}{\sig(1)<\cdots<\sig(j) \atop \sig(j+1)<\ldots<\sig(n)} \sig
	}
	
The key identities are as follows.

\begin{proposition}\label{2.2_Prop1}
We have
	\envHLineCFNmePd
	{\label{2.2_Prop1_EqRefIKZ}
		\lfgUT + \mo^n \lfgTN 
	}
	{
		\SmT{j=1}{n-1} \mo^{n-j-1} \lfgSHi{j} \lfgTO{j}
	}
	{\label{2.2_Prop1_EqCrucial}
		\lfgUT - \lfgEP \lfgT \lfgP \lfgT \lfgPi
	}
	{
		\lfgSHi{1} \bkR{ \lfgUT + \mo^n \lfgEP \lfgP \lfgT \lfgPi \lfgT \lfgP \lfgT \lfgPi }
		\lnAH[]
		-
		\mo^n \lfgEP \lfgP \bkR{ \lfgUT + \mo^n \lfgTN } \lfgPi \lfgSHi{1} \lfgP \lfgT \lfgPi \lfgT \lfgP \lfgT \lfgPi
		\nonumber
	}
\end{proposition}

Identity \refEq{2.2_Prop1_EqRefIKZ} was proved in \cite[\refLemA{2}]{IKZ06}.
Identity \refEq{2.2_Prop1_EqCrucial} was not proved in \cite{IKZ06},
	but we will show that it can be derived from the identity \cite[(8.5)]{IKZ06} in $\setZ[\gpSym{n+1}]$,
	\envHLineCm[2.2_PL_EqCyclicIKZ]
	{
		\lfgUTO{n+1} + \lfgSHiT{1}{n} \lfgCYCO{n+1}
	}
	{
		\lfgCYCO{n+1} \bkR[n]{ \lfgUTO{n+1} + \lfgSHiT{1}{n} \lfgTUN }
	}
	where $\lfgTUN=\lfgTUO{n+1}$ is the transposition $(1(1+n))$, 
	$\lfgCYCO{n+1}$ is the cyclic permutation of order $n+1$ defined by
	\envOTLineCmDef[2.2_PL_DefCycPerm]
	{
		\lfgCYCO{n+1}
	}
	{
		(12\ldots(n+1))			
	}
	{
		\left(\begin{array}{ccc|c}  0&   &   & 1  \\ \hline 1&  0 &   &   \\  &  \raMTx{-1}{\roMTx{-10}{\ddots}}  &  \!\raMTx{1}{\roMTx{-10}{\ddots}} &   \\  &   & \!1  & 0 \end{array}\right)
	\lnP{\in}
		\gpSym{n+1}
	\lnP{\subset}
		\gpGL{n+1}{\setZ}
	}
	and $\lfgSHiT{1}{n} \in \gpSym{n}$ is embedded into $\gpSym{n+1}$ 
	in a natural way,
	that is,
	a permutation $\SMatT{1&\ldots&n}{j_1&\ldots&j_n}$ of $\gpSym{n}$ 
	is identified with the permutation $\SMatT{1&\ldots&n&n+1}{j_1&\ldots&j_n&n+1}$ of $\gpSym{n+1}$ which fixes $n+1$.

To prove \refEq{2.2_Prop1_EqCrucial},
	we prepare notation for right actions of $\setZ[\gpGL{n}{\setZ}]$ and $\setZ[\gpSym{n+1}]$ 
	on the polynomial ring $\setR[x_1,\ldots,x_n]$ in $n$ variables with real coefficients by referring to \cite[\refSectA{8}]{IKZ06}.
The action of $S=\SmN a_j S_j\in\setZ[\gpGL{n}{\setZ}]$ on a polynomial $\Fc{f}{x_1,\ldots,x_n}$ is defined by
	\envHLineDefPd[2.2_PL_DefRAction]
	{
		\racFA{f}{S}{x_1,\ldots,x_n}
	}
	{
		\SmN a_j \Fc{f}{(x_1,\ldots,x_n)\cdot S_j^{-1}}
	}
Note that 
	\envM{
		\racFA{f}{\sig}{x_1,\ldots,x_n}
	}{
		\Fc{f}{x_{\sig^{-1}(1)},\ldots,x_{\sig^{-1}(n)}}
	}
	if $\sig\in\gpSym{n}$,
	and that this is a right action, i.e.,
	\envM{
		\racF{f}{(RS)}
	}{
		\racF{\racFrr{f}{R}}{S}
	}
	for all elements $R,S\in\setZ[\gpGL{n}{\setZ}]$.
The action of $\sig\in\gpSym{n+1}$ on $\setR[x_1,\ldots,x_n]$ is defined through the isomorphism
	\envHLinePt[2.2_PL_DefRActionExt]{\simeq}
	{
		\vphi : \setR[x_1,\ldots,x_n]
	}
	{
		\setR[y_1,\ldots,y_n,y_{n+1}]/(y_1+\cdots+y_n+y_{n+1}) 
	}
	given by
	\envHLinePd
	{
		\Fc{\vphi}{x_i}
	}
	{
		y_i
		\quad
		(1 \leq i \leq n)
	}
For example,
	the action $\sig\in\gpSym{n+1}$ on  $\Fc{p_1}{x_1,\ldots,x_n}=x_1\in \setR[x_1,\ldots,x_n]$ is
	\envHLine
	{
		\racFA{p_1}{\sig}{x_1,\ldots,x_n}
	}
	{
		\envCaseTPd{
			x_{\sig^{-1}(1)}
			&
			(1\leq \sig^{-1}(1) \leq n)
		}{
			-(x_1+\cdots+x_n)
			&
			(\sig^{-1}(1)=n+1)
		}
	}

We now prove \refEq{2.2_Prop1_EqCrucial}.
\envProof[identity \refEq{2.2_Prop1_EqCrucial}]{
By the definition of the action of $\gpSym{n+1}$,
	\envHLineCFCm
	{
		\racFA{p}{\lfgTUN}{x_1,\ldots,x_n}
	}
	{
		\Fc{p}{-(x_1+\cdots+x_n),x_2,\ldots,x_n}
	}
	{
		\racFA{p}{\lfgCYCO{n+1}}{x_1,\ldots,x_n}
	}
	{
		\Fc{p}{-(x_1+\cdots+x_n),x_1,\ldots,x_{n-1}}
	}
	where $\Fc{p}{x_1,\ldots,x_n}$ is a polynomial in $\setR[x_1,\ldots,x_n]$.
Since the map $\vphi$ is isomorphic,
	there exist matrices $\lfgTUOO{n}$ and $\lfgCYCT{n+1}{n}$ in $\gpGL{n}{\setZ}$ identified with $\lfgTUN$ and $\lfgCYCO{n+1}$ in $\gpSym{n+1}$,
	which are given by
	\envHLineCFLaaCm[2.2_Prop1Proof_DefTauCycProj1]
	{
		\lfgTUOO{n}
	}
	{
		\MatF{-1&& &}{ -1&1&&}{ \vdots& &\!\raMTx{2}{\roMTx{-9}{\ddots}}&  }{-1& & &1}
	}
	{
		\lfgCYCT{n+1}{n}
	}
	{
		\MatF{0&&&-1}{1&\roMTx{-9}{\ddots}&&\vdots}{&\!\raMTx{-1}{\roMTx{-9}{\ddots}}&0&-1}{&&1&-1}
	}
	respectively.
Note that the inverse matrices of $\lfgTUOO{n}$ and $\lfgCYCT{n+1}{n}$ are 
	\envHLineCFLaaPd[2.2_Prop1Proof_DefTauCycProj2]
	{
		\mopI[b]{\lfgTUOO{n}}
	}
	{
		\lfgTUOO{n}
	}
	{
		\mopI[b]{\lfgCYCT{n+1}{n}}
	}
	{
		\MatF{-1&1&&}{-1&0&\roMTx{-10}{\ddots}&}{\vdots&&\raMTx{1}{\roMTx{-10}{\ddots}}&1}{-1&&&0}
	\lnP{=}
		\lfgT \lfgCYCT{n+1}{n} \lfgT
	}
Using  \refEq{2.2_PL_DefP} and \refEq{2.2_PL_DefCycPerm},
	we can show that
	\envHLineCFLqqPd[2.2_Prop1Proof_EqTauCycProj]
	{
		\lfgTUOO{n}
	}
	{
		\lfgEP \lfgCYCO{n} \lfgPi \lfgT \lfgP \lfgT
	}
	{
		\lfgCYCT{n+1}{n}
	}
	{
		\lfgEP \lfgT\lfgPi \lfgT \lfgP 
	}
In fact we have by \refEq{2.2_PL_DefP} 
	\envOFLinePd
	{
		\lfgPi \lfgT \lfgP \lfgT
	}
	{
		\lfgPi (\lfgT \lfgP \lfgT)
	}
	{
		\MatF{1&& &}{ 1&\,1&&}{ \vdots&\vdots&\roMTx{-7}{\ddots}&  }{1&\,1& \cdots&1}
		\MatF{1&-1& &}{ &\ \roMTx{-0}{\ddots}&\roMTx{-0}{\ddots}&}{ & &1& -1}{& & &1}
	}
	{
		\MatF{1&-1& &}{ 1&0&\raMTx{2}{\roMTx{-12}{\ddots}}&}{ \vdots& &\raMTx{2}{\roMTx{-12}{\ddots}}&-1  }{1& & &0}
	}
We thus see that
	\envOFLineCm
	{
		\lfgCYCO{n} \lfgPi \lfgT \lfgP \lfgT
	}
	{
		\left(\begin{array}{cccc}  0&   &   & 1  \\ 1&  0 &   &   \\  &  \raMTx{-1}{\roMTx{-10}{\ddots}}  &  \!\raMTx{1}{\roMTx{-10}{\ddots}} &   \\  &   & \!1  & 0 \end{array}\right)
		\MatF{1&-1& &}{ 1&0&\raMTx{2}{\roMTx{-12}{\ddots}}&}{ \vdots& &\raMTx{2}{\roMTx{-12}{\ddots}}&-1  }{1& & &0}
	}
	{
		\MatF{1&& &}{ 1&-1&&}{ \vdots& &\!\raMTx{2}{\roMTx{-9}{\ddots}}&  }{1& & &-1}
	}
	{
		\lfgEP \lfgTUOO{n}
	}
	and 
	\envOFLinePd
	{
		\lfgT\lfgPi \lfgT \lfgP 
	}
	{
		\mopI[n]{ \lfgPi \lfgT \lfgP \lfgT }
	}
	{
		\MatF{1&-1& &}{ 1&0&\raMTx{2}{\roMTx{-12}{\ddots}}&}{ \vdots& &\raMTx{2}{\roMTx{-12}{\ddots}}&-1  }{1& & &0}^{-1}
	}
	{
		\MatF{0&&&1}{-1&\roMTx{-9}{\ddots}&&\vdots}{&\!\raMTx{-1}{\roMTx{-9}{\ddots}}&0&1}{&&-1&1}
	\lnP{=}
		\lfgEP \lfgCYCT{n+1}{n}
	}
These equations prove \refEq{2.2_Prop1Proof_EqTauCycProj}.

As is clear from the definition \refEq{2.2_PL_DefShuffle},
	the $1$-th and $(n-1)$-th shuffle elements can be rewritten in terms of cyclic permutations as
	\envHLineCFLaaCm[2.2_Prop1Proof_EqShuffleFirst]
	{
		\lfgSHiT{1}{n}
	}
	{
		\SmT[l]{j=1}{n} (j(j-1)\ldots1)
	}
	{
		\lfgSHiT{n-1}{n}
	}
	{
		\SmT[l]{j=1}{n} (j(j+1)\ldots n)
	}
	respectively.
We thus have
	\envMCm
	{
		\lfgSHiT{1}{n}\lfgCYCO{n}
	}
	{
		\lfgSHiT{n-1}{n}
	}
	and so obtain the following identity in $\setZ[\gpGL{n}{\setZ}]$
	by combining  \refEq{2.2_PL_EqCyclicIKZ} and \refEq{2.2_Prop1Proof_EqTauCycProj} under the identifications given by \refEq{2.2_PL_DefRActionExt}: 
	\envHLinePd[2.2_Prop1Proof_EqKeyIdentityIKZ1]
	{
		\lfgUTN + \lfgEP \lfgSHi{1} \lfgT\lfgPi \lfgT \lfgP 
	}
	{
		\lfgEP \lfgT\lfgPi \lfgT \lfgP  \bkR[n]{ \lfgUTN + \lfgEP \lfgSHi{n-1}  \lfgPi \lfgT \lfgP \lfgT }
	}
Equation \refEq{2.2_Prop1Proof_EqKeyIdentityIKZ1} can be restated by multiplying $\lfgEP \lfgPi \lfgT \lfgP \lfgT$ from the right as
	\envHLinePd[2.2_Prop1Proof_EqKeyIdentityIKZ2]
	{
		\lfgUTN - \lfgEP \lfgPi \lfgT \lfgP \lfgT
	}
	{
		\lfgSHi{1} - \lfgEP \lfgT\lfgPi \lfgT \lfgP \lfgSHi{n-1}  \lfgPi \lfgT \lfgP \lfgT \lfgPi \lfgT \lfgP \lfgT
	}

We define involutions $\pivIN$ and $\pivTRN$ of $\setZ[\gpGL{n}{\setZ}]$ 
	by the inversion and transpose of a matrix, 
	respectively,
	that is,
	for an element $S=\SmN a_j S_j\in\setZ[\gpGL{n}{\setZ}]$,
	we define the involutions as
	\envHLineCFLaaPd[2.2_Prop1Proof_DefInvol]
	{
		\pivIO{S}
	}
	{
		\SmN a_j \mopI{S_j} 
	}
	{
		\pivTRO{S}
	}
	{
		\SmN a_j \mopTP{S_j} 	
	}
Note that $\pivIO{\sig}=\pivTRO{\sig}$ if $\sig\in\gpSym{n}$,
	in particular,
	$ \Fc{\pivTRN\circ\pivIN}{\lfgSHi{j}} = \lfgSHi{j}$.
Since 
	$ \Fc{\pivTRN}{P^{\pm}} = \lfgT P^{\pm} \lfgT$
	and
	$\lfgT\lfgSHi{n-1}\lfgT = \lfgSHi{1}$,
	we see that
	\envOTLineThCm
	{
		\Fc{ \pivTRN\circ\pivIN }{ \lfgUTN - \lfgEP \lfgPi \lfgT \lfgP \lfgT }
	}
	{
		\Fc{ \pivTRN }{ \lfgUTN - \lfgEP \lfgT \lfgPi \lfgT \lfgP  }
	}
	{
		\lfgUT - \lfgEP \lfgT \lfgP \lfgT \lfgPi
	}
	and
	\envLineFPd
	{
		\Fc{ \pivTRN\circ\pivIN }{ \lfgSHi{1} - \lfgEP \lfgT\lfgPi \lfgT \lfgP \lfgSHi{n-1}  \lfgPi \lfgT \lfgP \lfgT \lfgPi \lfgT \lfgP \lfgT }
	}
	{
		\lfgSHi{1}  - \Fc{ \pivTRN }{ \lfgEP \lfgT\lfgPi \lfgT \lfgP \lfgT\lfgPi \lfgT \lfgP  \Fc{\pivIN}{ \lfgSHi{n-1} } \lfgPi \lfgT \lfgP \lfgT }
	}
	{
		\lfgSHi{1}  - \lfgEP \lfgT (\lfgT\lfgP\lfgT) \lfgT (\lfgT\lfgPi\lfgT) \lfgSHi{n-1} (\lfgT\lfgP\lfgT) \lfgT (\lfgT\lfgPi\lfgT) \lfgT (\lfgT\lfgP\lfgT) \lfgT (\lfgT\lfgPi\lfgT) \lfgT
	}
	{
		\lfgSHi{1} - \lfgEP \lfgP \lfgT \lfgPi \lfgSHi{1} \lfgP \lfgT \lfgPi \lfgT \lfgP \lfgT \lfgPi
	}
Thus,
	by applying the involuation $\pivTRN\circ\pivIN$ to both sides of \refEq{2.2_Prop1Proof_EqKeyIdentityIKZ2},
	we obtain
	\envHLinePd[2.2_Prop1Proof_EqKeyIdentityOrigin]
	{
		\lfgUT - \lfgEP \lfgT \lfgP \lfgT \lfgPi
	}
	{
		\lfgSHi{1} - \lfgEP \lfgP \lfgT \lfgPi \lfgSHi{1} \lfgP \lfgT \lfgPi \lfgT \lfgP \lfgT \lfgPi
	}
The second term of the right-hand side of \refEq{2.2_Prop1Proof_EqKeyIdentityOrigin} is calculated as
	\envLineThCm
	{\quad
		\lfgEP \lfgP \lfgT \lfgPi \lfgSHi{1} \lfgP \lfgT \lfgPi \lfgT \lfgP \lfgT \lfgPi
	}
	{
		\mo^n \lfgEP \lfgP \bkB{ \bkR{ \lfgUT + \mo^n \lfgTN } - \lfgUT } \lfgPi \lfgSHi{1} \lfgP \lfgT \lfgPi \lfgT \lfgP \lfgT \lfgPi
	}
	{
		-
		\mo^n \lfgEP \lfgSHi{1} \lfgP \lfgT \lfgPi \lfgT \lfgP \lfgT \lfgPi 
		+
		\mo^n \lfgEP \lfgP \bkR{ \lfgUT + \mo^n \lfgTN } \lfgPi \lfgSHi{1} \lfgP \lfgT \lfgPi \lfgT \lfgP \lfgT \lfgPi
	}
	and so we obtain \refEq{2.2_Prop1_EqCrucial}.
}

\section{Proofs} \label{sectThree}
We denote by $K[x_1,\ldots,x_n]_{(m)}$ the set of homogeneous polynomials of degree $m$ in $n$ variables with coefficients in a field $K$ of characteristic zero. 
We can consider $K[x_1,\ldots,x_n]_{(m)}$ as the vector space over $\setQ$,
	\envHLineCm[3_PL_DefModule]
	{
		K[x_1,\ldots,x_n]_{(m)} 
	}
	{
		\tpOPs{k_1,\ldots,k_n\geq0}{k_1+\cdots+k_n=m} K x_1^{k_1} \cdots x_n^{k_n} 
	}
	and so we also denote by $V[x_1,\ldots,x_n]_{(m)}$ the vector space over $\setQ$ given by \refEq{3_PL_DefModule} with $K=V$
	for any space $V$.

We define the generating functions $\ggfcZHbb=\ggfcZH[]{n}{}$ and $\ggfcZSbb=\ggfcZS[]{n}{}$ of RMZVs of harmonic and shuffle types by
	\envHLineDef[3_PL_DefGGgfcRMZVh]
	{
		\ggfcZH{n}{x_1,\ldots,x_n}
	}
	{
		\Sm{l_1,\ldots,l_n\geq1} \fcZH{l_1,\ldots,l_n} x_1^{l_1-1} \cdots x_n^{l_n-1}
	}
	and
	\envHLineCmDef[3_PL_DefGGgfcRMZVs]
	{
		\ggfcZS{n}{x_1,\ldots,x_n}
	}
	{
		\Sm{l_1,\ldots,l_n\geq1} \fcZS{l_1,\ldots,l_n} x_1^{l_1-1} \cdots x_n^{l_n-1}
	}
	respectively.
Because our interests are RMZVs of weight $l$,
	we also define their generating functions $\gfcZHbb{l}=\gfcZH[]{l}{n}{}$ and $\gfcZSbb{l}=\gfcZS[]{l}{n}{}$ by
	\envHLineDef[3_PL_DefGgfcRMZVh]
	{
		\gfcZH{l}{n}{x_1,\ldots,x_n}
	}
	{
		\tpSm{l_1,\ldots,l_n\geq1}{l_1+\cdots+l_n=l} \fcZH{l_1,\ldots,l_n} x_1^{l_1-1} \cdots x_n^{l_n-1}
	}
	and
	\envHLineCmDef[3_PL_DefGgfcRMZVs]
	{
		\gfcZS{l}{n}{x_1,\ldots,x_n}
	}
	{
		\tpSm{l_1,\ldots,l_n\geq1}{l_1+\cdots+l_n=l} \fcZS{l_1,\ldots,l_n} x_1^{l_1-1} \cdots x_n^{l_n-1}
	}
	which are the homogeneous parts of degree $l-n$ of $\ggfcZH[]{n}{}$ and $\ggfcZS[]{n}{}$, respectively.
In a similar fashion to usage of $\fcZR[]{}$, $\fcZH[]{}$, and $\fcZS[]{}$,
	we use the same symbol $\gfcZRbb{l}$ for both $\gfcZHbb{l}$ and $\gfcZSbb{l}$ 
	when we consider a congruence relation modulo a vector space including $\wspPZV{l}[x_1,\ldots,x_n]_{(l-n)}$.
We omit the subscript $(l-n)$ for brevity (e.g., we use $\wspPZV{l}[x_1,\ldots,x_n]$ instead of $\wspPZV{l}[x_1,\ldots,x_n]_{(l-n)}$) 
	when the homogeneous degree in question can be inferred from the weight and depth.

We need a lemma to prove \refThm[s]{1_Thm1} and \ref{1_Thm2}.
\begin{lemma}\label{3_Lem1}
We have
	\envHLineCFNmePdPt{\equiv}
	{\label{3_Lem1_EqShu}
		\racF{\gfcZRbb{l}}{\lfgP(\lfgUT+\mo^n\lfgT)}
	}
	{
		0
		\qquad
		\modd \wspPZV{l}[x_1,\ldots,x_n]
	}
	{\label{3_Lem1_EqHar}
		\racF{\gfcZRbb{l}}{(\lfgUT - \lfgEP \lfgT \lfgP \lfgT \lfgPi)}
	}
	{
		0
		\qquad
		\modd (\wspZVd{l}{n-1}+\wspPZV{l})[x_1,\ldots,x_n]
	}
\end{lemma}
\envProof{
We may prove
	\envHLineCFNmePdPt{\equiv}
	{\label{3_Lem1Proof_EqShu}
		\racF{\gfcZRbb{l}}{\lfgP\lfgSHi{j}}
	}
	{
		0
		\quad
		\modd \wspPZV{l}[x_1,\ldots,x_n]
		\qquad
		(1\leq j \leq n-1)
	}
	{\label{3_Lem1Proof_EqHar}
		\racF{\gfcZRbb{l}}{\lfgSHi{1}}
	}
	{
		0
		\quad
		\modd (\wspZVd{l}{n-1}+\wspPZV{l})[x_1,\ldots,x_n]
	}
In fact,
	\refEq{3_Lem1_EqShu} follows from \refEq{2.2_Prop1_EqRefIKZ} and \refEq{3_Lem1Proof_EqShu},
	and
	\refEq{3_Lem1_EqHar} follows from \refEq{2.2_Prop1_EqCrucial}, \refEq{3_Lem1_EqShu}, and \refEq{3_Lem1Proof_EqHar}.
(Note that \refEq{3_Lem1Proof_EqShu} and \refEq{3_Lem1Proof_EqHar} 
	are essentially shown in the proofs of the assertions (ii) and (i) in the proof of \cite[\refThmA{6}]{IKZ06}, 
	respectively.)

We first prove \refEq{3_Lem1Proof_EqShu}.
Let $\lgP{l_1,\ldots,l_n}{t}$ be the multiple polylogarithm defined by
	\envHLineDefPd
	{
		\lgP{l_1,\ldots,l_n}{t}
	}
	{
		\Sm{m_1>\cdots>m_n>0} \opF{t^{m_1}}{ \pw{m_1}{l_1}\cdots\pw{m_n}{l_n} }
	}
We define $\idFc{\zeta}{\ep}{l_1,\ldots,l_n}=\lgP{l_1,\ldots,l_n}{1-\ep}$
	and its generating function as 
	\envHLineDefPd
	{
		\pggfcZe{n}{x_1,\ldots,x_n}
	}
	{
		\Sm{l_1,\ldots,l_n\geq1} \idFc{\zeta}{\ep}{l_1,\ldots,l_n} x_1^{l_1-1} \cdots x_n^{l_n-1}
	}
Let $\popSP{f}$ denote $\racF{f}{\lfgPO{n}}$ for a function $\Fc{f}{x_1,\ldots,x_n}$ in $n$ variables.
In the proof of \cite[\refThmA{6}]{IKZ06},
	\etalTx{Ihara} showed that
	\envHLinePd[3_Lem1Proof_Eq1Shu]
	{
		\racFA{\popSPr{\pggfcZe[]{n}{}}}{\lfgSHi{j}}{x_1,\ldots,x_n}
	}
	{
		\popSPrT{\pggfcZe[]{j}{}}{x_1,\ldots,x_j} \popSPrT{\pggfcZe[]{n-j}{}}{x_{j+1},\ldots,x_n}
	}
Since
	the value $\fcZS{l_1,\ldots,l_n}$ is equal to 
	the constant term of the asymptotic expansion of the multiple polylogarithm $\lgP{l_1,\ldots,l_n}{t}$ as $t\nearrow1$ (see \cite[\refSecA{2}]{IKZ06}),
	we obtain by comparing the constant terms as $\ep\searrow0$ on both sides of \refEq{3_Lem1Proof_Eq1Shu} that 
	\envHLinePd[3_Lem1Proof_Eq1Shu2]
	{
		\racFA{\popSPr{\ggfcZS[]{n}{}}}{\lfgSHi{j}}{x_1,\ldots,x_n}
	}
	{
		\popSPrT{\ggfcZS[]{j}{}}{x_1,\ldots,x_j} \popSPrT{\ggfcZS[]{n-j}{}}{x_{j+1},\ldots,x_n}
	}
The homogeneous part of degree $l-n$ of the left-hand side of \refEq{3_Lem1Proof_Eq1Shu2} 
	is equal to $\racF{\popSPr{\gfcZRbb{l}}}{\lfgSHi{j}}=\racF{\gfcZRbb{l}}{\lfgP\lfgSHi{j}}$,
	and 
	that of the right-hand side of \refEq{3_Lem1Proof_Eq1Shu2} is congruent to $0$ modulo $\wspPZV{l}[x_1,\ldots,x_n]$.
Therefore \refEq{3_Lem1Proof_Eq1Shu2} proves \refEq{3_Lem1Proof_EqShu}. 

We next prove \refEq{3_Lem1Proof_EqHar}.
It follows from \refEq{2_PL_DefHarProdMain} and induction on $n$ that 
	\envHLinePd
	{
		z_{l_1} * z_{l_2} \cdots z_{l_n}
	}
	{
		\SmT{j=1}{n} z_{l_2} \cdots z_{l_j} z_{l_1} z_{l_{j+1}} \cdots z_{l_n}
		+
		\SmT{j=2}{n} z_{l_2} \cdots z_{l_{j-1}} z_{l_1+l_j} z_{l_{j+1}} \cdots z_{l_n}
	}
Applying the map $\mpEVH[]{}$ to both sides of this equation and substituting $T=0$,
	we obtain 
	\envHLinePd[3_Lem1Proof_Eq1Har]
	{
		\fcZH{l_1}\fcZH{l_2,\ldots,l_n} 
	}
	{
		\SmT{j=1}{n}  \fcZH{l_2,\ldots,l_j, l_1 ,l_{j+1},\ldots,l_n} 
		\lnAH
		+
		\SmT{j=2}{n} \fcZH{l_2,\ldots,l_{j-1}, l_1+l_j ,l_{j+1},\ldots,l_n} 
	}
We easily see that
	the left-hand side is congruent to $0$ modulo $\wspPZV{l}$.
We also see that
	the first sum in the right-hand side of \refEq{3_Lem1Proof_Eq1Har} can be rewritten as
	\envHLine
	{
		\SmT[l]{j=1}{n}  \fcZH{l_2,\ldots,l_j, \osTx{l_1}{j\,th} ,l_{j+1},\ldots,l_n} 
	}
	{
		\racFA{\fcZH[]{}}{\lfgSHi{1}}{l_1,\ldots,l_n}
	}
	by using \refEq{2.2_Prop1Proof_EqShuffleFirst},
	and that
	the second sum is congruent to $0$ modulo $\wspZVd{l}{n-1} + \wspPZV{l}$ by \refEq{2.1_Prop1_Eq2}.
We thus derive from \refEq{2.1_Prop1_Eq1} and \refEq{3_Lem1Proof_Eq1Har} that 
	\envHLineCmPt{\equiv}
	{
		\racFA{\fcZR[]{}}{\lfgSHi{1}}{l_1,\ldots,l_n}
	} 
	{
		0
		\qquad
		\wspZVd{l}{n-1}+\wspPZV{l}
	}
	which proves \refEq{3_Lem1Proof_EqHar} by the definitions of the generating functions $\gfcZHbb{l}$ and $\gfcZSbb{l}$.
}

We are now able to prove \refThm{1_Thm1}.
\envProof[\refThm{1_Thm1}]{
A direct calculation shows that
	\envHLineCm
	{
		\lfgUT+\mo^n\lfgEP\lfgT
	}
	{
		\lfgP (\lfgUT+\mo^n\lfgT) \lfgPi 
		-
		\mo^n ( \lfgUT - \lfgEP \lfgT \lfgP \lfgT \lfgPi ) \lfgP \lfgT\lfgPi 
	}
	which together with \refEq{3_Lem1_EqShu} and \refEq{3_Lem1_EqHar} gives
	\envHLinePdPt{\equiv}
	{
		\gfcZRbb{l}
	}
	{
		\mo^{n-1} \racF{\gfcZRbb{l}}{\lfgEP\lfgT}
		\qquad
		\modd (\wspZVd{l}{n-1}+\wspPZV{l})[x_1,\ldots,x_n]
	}
Because
	\envMTh{
		\racFA{\gfcZRbb{l}}{\lfgEP\lfgT}{x_1,\ldots,x_n}
	}{
		\racFA{\gfcZRbb{l}}{\lfgT}{-x_1,\ldots,-x_n}
	}{
		\mo^{l-n} \racFA{\gfcZRbb{l}}{\lfgT}{x_1,\ldots,x_n}
	}
	by the homogeneity of $\gfcZRbb{l}$,
	we obtain
	\envHLinePdPt{\equiv}
	{
		\gfcZRbb{l}
	}
	{
		\mo^{l-1} \racF{\gfcZRbb{l}}{\lfgT}
		\qquad
		\modd (\wspZVd{l}{n-1}+\wspPZV{l})[x_1,\ldots,x_n]
	}
Comparing the coefficient of $x_1^{l_1-1}\ldots x_n^{l_n-1}$ on both sides of this equation,
	we deduce that
	\envMCmPt{\equiv}
	{
		\fcZR{\lbfL}
	}
	{
		\mo^{l-1} \fcZR{\lbfL\cdot\lfgT}
		\ \modd \wspZVd{l}{n-1}+\wspPZV{l}
	}
	where $\lbfL=(l_1,\ldots,l_n)$.
Since $\lbfL\cdot\lfgT=\popR{\lbfL}$ by the definition of $\vee$,
	this congruence identity can be rewritten as 
	\envHLineCmPt[3_Thm1Proof_EqRefRslt]{\equiv}
	{
		\fcZR{\lbfL}
	}
	{
		\mo^{l-1} \fcZR{\popR{\lbfL}}
		\qquad
		\modd \wspZVd{l}{n-1}+\wspPZV{l}
	}
	which completes the proof of the first claim of \refThm{1_Thm1}.
	
The claim left immediately follows from the fact that 
	equating \refEq{1_PL_EqqRefRsltFromIKZ2} and \refEq{3_Thm1Proof_EqRefRslt} proves \refEq{1_PL_EqqParityResult}
	because $\wspZVd{l}{n-1}\subset\wspZVd{l}{<n}$ (see also \refRem{3_Rem1} below).
}
\begin{remark}\label{3_Rem1}
We will briefly prove \refEq{1_PL_EqqRefRsltFromIKZ1} by using \refEq{2.2_Prop1Proof_EqKeyIdentityIKZ1} instead of \refEq{2.2_Prop1Proof_EqKeyIdentityOrigin}
	in the same manner as the proof of \refEq{3_Thm1Proof_EqRefRslt}.	
Since \refEq{2.2_Prop1Proof_EqKeyIdentityIKZ1} is equivalent to \refEq{2.2_PL_EqCyclicIKZ} or \cite[(8.5)]{IKZ06},
	the proof of \refEq{1_PL_EqqRefRsltFromIKZ1} yields a restatement of the proof of the parity result \refEq{1_PL_EqqParityResult} 
	given by \etalTx{Ihara} \cite{IKZ06} in the view point of our strategy.

Multiplying both sides of \refEq{2.2_Prop1Proof_EqKeyIdentityIKZ1} by $\lfgEP\lfgP\lfgT$ from the left and using $\lfgT\lfgSHi{1}\lfgT = \lfgSHi{n-1}$,
	we can obtain
	\envHLinePd[3_Rem1_EqKeyIdentityIKZ1]
	{
		\lfgEP \lfgP \lfgT - \lfgT \lfgP
	}
	{
		-
		\lfgP \lfgSHi{n-1} \lfgPi \lfgT \lfgP 
		+
		\lfgEP \lfgT \lfgP \lfgSHi{n-1}  \lfgPi \lfgT \lfgP \lfgT 
	}
We see from $\lfgT=\mo^n \bkB{ (\lfgUT + \mo^n\lfgT) - \lfgUT }$ that
	\envHLineCFCm
	{
		\lfgEP \lfgP \lfgT - \lfgT \lfgP
	}
	{
		\mo^n \lfgEP \bkB{ \lfgP (\lfgUT + \mo^n\lfgT) - (\lfgUT + \mo^n\lfgEP\lfgT) \lfgP }
	}
	{
		\lfgEP \lfgT \lfgP \lfgSHi{n-1}  \lfgPi \lfgT \lfgP \lfgT 
	}
	{
		-
		\mo^n \lfgEP \lfgP \lfgSHi{n-1} \lfgPi \lfgT \lfgP \lfgT  
		\lnAH
		+ 
		\mo^n \lfgEP (\lfgUT + \mo^n\lfgT) \lfgP \lfgSHi{n-1}  \lfgPi \lfgT \lfgP \lfgT 
	}
	which together with \refEq{3_Rem1_EqKeyIdentityIKZ1} give
	\envHLinePd[3_Rem1_EqKeyIdentityIKZ2]
	{
		(\lfgUT + \mo^n\lfgEP\lfgT)\lfgP
	}
	{
		\mo^n \lfgEP \lfgP \lfgSHi{n-1} \lfgPi \lfgT \lfgP (\lfgUT + \mo^n\lfgEP\lfgT)
		\lnAH[]
		-
		(\lfgUT + \mo^n\lfgT) \lfgP \lfgSHi{n-1}  \lfgPi \lfgT \lfgP \lfgT 
		+
		\lfgP (\lfgUT + \mo^n\lfgT)
	\nonumber
	}
We know from \refEq{1_PL_EqqRefRsltFromIKZ2} (or \cite[(8.6)]{IKZ06}) that
	\envHLinePdPt[3_Rem1_EqReflectHar]{\equiv}
	{
		\racF{\gfcZRbb{l}}{(\lfgUT+\mo^n\lfgT)}
	}
	{
		0
		\qquad 
		\modd (\wspZVd{l}{<n}+\wspPZV{l})[x_1,\ldots,x_n]
	}
Thus we can derive 
	from \refEq{3_Rem1_EqKeyIdentityIKZ2} together with \refEq{3_Lem1_EqShu}, \refEq{3_Lem1Proof_EqShu}, and \refEq{3_Rem1_EqReflectHar} that
	\envHLineCmPt{\equiv}
	{
		\racF{\gfcZRbb{l}}{(\lfgUT+\mo^n\lfgEP\lfgT)\lfgP}
	}
	{
		0
		\qquad 
		\modd (\wspZVd{l}{<n}+\wspPZV{l})[x_1,\ldots,x_n]
	}
	or equivalently,
	\envHLineCmPt{\equiv}
	{
		\gfcZRbb{l}
	}
	{
		\mo^{n-1} \racF{\gfcZRbb{l}}{\lfgEP\lfgT}
		\qquad
		\modd (\wspZVd{l}{<n}+\wspPZV{l})[x_1,\ldots,x_n]
	}
	from which we obtain \refEq{1_PL_EqqRefRsltFromIKZ1} in the same way as the proof of \refThm{1_Thm1}.
\end{remark}
We now prove \refThm{1_Thm2}.

\envProof[\refThm{1_Thm2}]{
For any congruence equation $\equiv$ in this proof,
	we suppose the vector space giving the equivalent relation to be $(\wspZVd{l}{n-1}+\wspPZV{l})[x_1,\ldots,x_n]$.

We see from $\mopTP{\lfgP^{\pm}} = \lfgT\lfgP^{\pm}\lfgT$ and \refEq{2.2_Prop1Proof_EqTauCycProj} that
	\envOTLineCm[3_Thm2Proof_EqTPTPi1]
	{
		\lfgEP \lfgT \lfgP \lfgT \lfgPi
	}
	{
		\mopTP[n]{\lfgEP \lfgT\lfgPi \lfgT \lfgP }
	}
	{
		\mopTP[m]{\lfgCYCT{n+1}{n}}
	}
	which together with \refEq{2.2_Prop1Proof_DefTauCycProj2} gives
	\envOFLinePd[3_Thm2Proof_EqTPTPi2]
	{
		\mopI[n]{\lfgEP \lfgT \lfgP \lfgT \lfgPi}
	}
	{
		\mopI[b]{\mopTP{\lfgCYCT{n+1}{n}} }
	}
	{
		\mopTP[b]{\mopI[b]{\lfgCYCT{n+1}{n}}}
	}
	{
		\MatF{-1&-1&-1&-1}{1&0&&}{&\raMTx{-2}{\roMTx{-6}{\ddots}}&\raMTx{0}{\roMTx{-6}{\ddots}}&}{&&1&0}
	}
We obtain from \refEq{3_Lem1_EqHar} and \refEq{3_Thm2Proof_EqTPTPi2} that 
	\envOTLinePdPt{\equiv}
	{
		\gfcZRb{l}{x_1,\ldots,x_n}
	}
	{
		\gfcZRb{l}{(x_1,\ldots,x_n) \mopI[n]{\lfgEP \lfgT \lfgP \lfgT \lfgPi} }
	}
	{
		\gfcZRb{l}{x_2-x_1,\ldots,x_n-x_1, -x_1}
	}
Replacing $x_i$ by $x_i-x_{n+1}$ for every integer $i$ with $1\leq i \leq n$,
	we can rewrite this as 
	\envHLineThCmPt{\equiv}
	{
		\gfcZRb{l}{x_1-x_{n+1},\ldots,x_n-x_{n+1}}
	}
	{
		\gfcZRb{l}{x_2-x_1,\ldots,x_{n+1}-x_1}
	}
	{
		\gfcZRb{l}{x_{c(1)}-x_{c(n+1)},\ldots,x_{c(n)}-x_{c(n+1)}}
	}
	where $c$ means the cyclic permutation $(12\ldots(n+1))$ in $\gpSym{n+1}$.
We thus have
	\envHLinePt[3_Thm2Proof_EqGfcCyclic]{\equiv}
	{\hspace{-10pt}
		\gfcZRb{l}{x_1-x_{n+1},\ldots,x_n-x_{n+1}}
	}
	{
		\gfcZRb{l}{x_{c^i(1)}-x_{c^i(n+1)},\ldots,x_{c^i(n)}-x_{c^i(n+1)}}
	}
	for each integer $i$ with $1\leq i \leq n$.
Identity \refEq{3_Thm2Proof_EqGfcCyclic} with $x_i=x_{n+1}=0$ yields
	\envHLineCmPt{\equiv}
	{
		\gfcZRb{l}{x_1,\ldots,x_{i-1},\osTx{0}{i\,th},x_{i+1},\ldots,x_n}
	}
	{
		\gfcZRb{l}{x_{i+1},\ldots,x_n,\osTx{0}{j\,th},x_1,\ldots,x_{i-1}}
	}
	where $j=n+1-i$.
By comparing the coefficient of $x_1^{l_1-1}\cdots x_n^{l_n-1}$ on both sides of this equation,
	we obtain \refEq{1_Thm2_EqqRefRslt1}.
	
We see from \refThm{1_Thm1} that 
	\envOTLineCmPt{\equiv}
	{
		\fcZR{\pbfL_i,1,\pbfK_i}
	}
	{
		\mo^{l-1} \fcZR{\popR[n]{\pbfL_i,1,\pbfK_i}}
	}
	{
		\mo^{l-1} \fcZR{ \popR{\pbfK_i},1,\popR{\pbfL_i}}
	}
	which together with \refEq{1_Thm2_EqqRefRslt1} proves \refEq{1_Thm2_EqqRefRslt2}.
}
Let $\bkB{1}^{m}$ be $m$ repetitions of $1$,
	where $\bkB{1}^{0}$ means the empty set. 
We prepare a lemma to prove \refCor{1_Cor1}.
\begin{lemma}\label{3_Lem2}
Let $l,p,q$ be integers with $l\geq p+q$ and $p,q\geq1$.
If weight $l$ is even, 
	\envHLinePdPt[3_Lem2_Eq]{\equiv}
	{
		\tpSm{l_1\geq2, l_2,\ldots,l_p\geq1}{l_1+\cdots+l_p=l-q+1} \fcZ{l_1,\ldots,l_p,\bkB{1}^{q-1}}
	}
	{
		0
		\qquad
		\modd \wspZVd{l}{p+q-2} + \wspPZV{l}
	}
\end{lemma}
\envProof{
We define a sum $\pnmSr{l,p,q}$ of RMZVs of weight $l$ by
	\envHLineDefPd[3_Lem2Proof_DefSumRepeatOne]
	{
		\pnmSr{l,p,q}
	}
	{
		\tpSm{l_1, \ldots,l_p\geq1}{l_1+\cdots+l_p=l-q+1} \fcZR{l_1,\ldots,l_p,\bkB{1}^{q-1}}
	}
If $q>1$,		
	we obtain by \refEq{1_Thm2_EqqRefRslt2} with $i=p+1$ 
	\envHLineFPdPt[3_Lem2Proof_EqSmRMZVhelp]{\equiv}
	{
		\pnmSr{l,p,q}
	}
	{
		\mo^{l-1} \tpSm{l_1, \ldots,l_p\geq1}{l_1+\cdots+l_p=l-q+1} \fcZR{ \popR[n]{l_1,\ldots,l_p},1, \popR[n]{ \usbTx{q-2}{1,\ldots,1} } }
	}
	{
		\mo^{l-1} \tpSm{l_1, \ldots,l_p\geq1}{l_1+\cdots+l_p=l-q+1} \fcZR{l_p,\ldots,l_1,\bkB{1}^{q-1}}
	}
	{
		\mo^{l-1} \pnmSr{l,p,q}
		\qquad
		\modd \wspZVd{l}{p+q-2} + \wspPZV{l}
	}
We can also calculate as \refEq{3_Lem2Proof_EqSmRMZVhelp} when $q=1$ by using \refThm{1_Thm1}.
Therefore it holds that
	\envHLinePdPt[3_Lem2Proof_EqSmRMZV]{\equiv}
	{
		\pnmSr{l,p,q}
	}
	{
		0
		\qquad
		\modd \wspZVd{l}{p+q-2} + \wspPZV{l}
		\qquad\ 
		(\text{$l$ is even})
	}

Let $\pnmS{l,p,q}$ be the left-hand side of \refEq{3_Lem2_Eq}.
We see from \refEq{1_Thm2_EqqRefRslt1} with $i=1$ that
	\envHLineThCmPt{\equiv}
	{
		\tpSm{l_2,\ldots,l_p\geq1}{l_2+\cdots+l_p=l-q} \fcZR{1,l_2,\ldots,l_p,\bkB{1}^{q-1}}
	}
	{
		\tpSm{l_2,\ldots,l_p\geq1}{l_2+\cdots+l_p=l-q} \fcZR{l_2,\ldots,l_p,\bkB{1}^q}
	}
	{
		\pnmSr{l,p-1,q+1}
		\qquad
		\modd \wspZVd{l}{p+q-2} + \wspPZV{l}
	}
	which  gives
	\envHLineThPdPt[3_Lem2Proof_EqSmMZV]{\equiv}
	{\hspace{-15pt}
		\pnmS{l,p,q}
	}
	{
		\tpSm{l_1, \ldots,l_p\geq1}{l_1+\cdots+l_p=l-q+1} \fcZR{l_1,\ldots,l_p,\bkB{1}^{q-1}} 
		- 
		\tpSm{l_2,\ldots,l_p\geq1}{l_2+\cdots+l_p=l-q} \fcZR{1,l_2,\ldots,l_p,\bkB{1}^{q-1}}
	}
	{
		\pnmSr{l,p,q} - \pnmSr{l,p-1,q+1}
		\qquad
		\modd \wspZVd{l}{p+q-2} + \wspPZV{l}
	}
Combining \refEq{3_Lem2Proof_EqSmRMZV} and \refEq{3_Lem2Proof_EqSmMZV} proves \refEq{3_Lem2_Eq}.
}

We give a proof of \refCor{1_Cor1}.

\envProof[\refCor{1_Cor1}]{
We define a sum $\pnmR{l,a,b}$ of MZVs of weight $l$ by
	\envHLineDef[3_Lem2Proof_DefSumFirstArgument]
	{
		\pnmR{l,a,b}
	}
	{
		\tpSm{l_1\geq a, l_2,\ldots,l_b\geq1}{l_1+\cdots+l_b=l} \fcZ{l_1,\ldots,l_b}
	}
	for integers $a,b$ with  $a\geq2$, $b\geq 1$, and $l\geq a+b-1$.
Let $\pnmS{l,p,q}$ be the left-hand side of \refEq{3_Lem2_Eq} as in the proof of \refLem{3_Lem2}.
The restricted sum formula \cite{ELO09} can be rewritten as 
	\envHLineCm[3_Cor1Proof_EqRSumFormula]
	{
		\pnmR{l,m+1,n}
	}
	{
		\pnmS{l,l-m-n+1,n}
	}
	where $m,n\geq1$ and $l\geq m+n$.
We recall the duality theorem (see \cite{Zagier94})
	\envHLine[3_Lem2Proof_EqDuality]
	{
		\fcZ{\pbfK}
	}
	{
		\fcZ{\pbfK'}
	}
	for any dual pair $(\pbfK,\pbfK')$ of index sets defined by
	\envM
	{
		\pbfK
	}
	{
		(a_1+1,\bkB{1}^{b_1-1}, \ldots, a_h+1,\bkB{1}^{b_h-1})
	}
	and
	\envMPd
	{
		\pbfK'
	}
	{
		(b_h+1,\bkB{1}^{a_h-1}, \ldots, b_1+1,\bkB{1}^{a_1-1})
	}
Ohno's theorem \cite{Ohno99} for the dual identity $\fcZ{m+1,\bkB{1}^{n-1}}=\fcZ{n+1,\bkB{1}^{m-1}}$ states that
	\envHLinePd[3_Cor1Proof_EqOhnoRel]
	{
		\pnmR{l,m+1,n}
	}
	{
		\pnmR{l,n+1,m}
	}
Equating \refEq{3_Cor1Proof_EqRSumFormula} and \refEq{3_Cor1Proof_EqOhnoRel} yields
	\envOTLinePd[3_Cor1Proof_EqRSumFormulaFromOnho]
	{
		\pnmR{l,m+1,n}
	}
	{
		\pnmS{l,l-m-n+1,m}
	}
	{
		\pnmS{l,l-m-n+1,n}
	}
Therefore it follows from \refEq{3_Lem2_Eq} and \refEq{3_Cor1Proof_EqRSumFormulaFromOnho} that
	\envHLineCmPt[3_Cor1Proof_EqqRestrictedOhnoRel]{\equiv}
	{
		\pnmR{l,m+1,n}
	}
	{
		0
		\qquad
		\modd \wspZVd{l}{l-d-1} + \wspPZV{l}
	}
	where $d\in\SetO{m,n}$.
We know from the duality theorem \refEq{3_Lem2Proof_EqDuality} that
	$\wspZVd{l}{l-d-1}=\wspZVd{l}{d+1}$,
	and so we obtain \refEq{1_Cor1_EqqRestrictedOhnoRel}.
}

\begin{table}[!t]\renewcommand{\arraystretch}{1.5}
\begin{center}{\scriptsize
	\begin{tabular}{|c|l|}  \hline
		\bfTx{mod} 					&\bfTx{Identity}																			\\\hline  
		$\wspZVd{3}{1} + \wspPZV{3}$     	&$\fcZ{2,1}\equiv \fcZR{1,2}$																\\\hline
		$\wspZVd{4}{1} + \wspPZV{4}$     	&$\fcZ{3,1}\equiv -\fcZ{3,1},\fcZR{1,3}$														\\\cdashline{1-2}
		$\wspZVd{4}{2} + \wspPZV{4}$     	&$\fcZ{2,1,1}\equiv -\fcZR{2,1,1},\pm\fcZR{1,2,1},\fcZR{1,1,2}$										\\\hline
		$\wspZVd{5}{1} + \wspPZV{5}$     	&$\fcZ{4,1}\equiv \fcZR{1,4}$																\\\cdashline{1-2}
		$\wspZVd{5}{2} + \wspPZV{5}$     	&$\fcZ{3,1,1}\equiv \fcZR{1,3,1},\fcZR{1,1,3}$													\\
								     	&$\fcZ{2,2,1}\equiv \fcZR{1,2,2}$															\\\cdashline{1-2}
		$\wspZVd{5}{3} + \wspPZV{5}$     	&$\fcZ{2,1,1,1}\equiv \fcZR{1,2,1,1},\fcZR{1,1,2,1},\fcZR{1,1,1,2}$									\\\hline
		$\wspZVd{6}{1} + \wspPZV{6}$     	&$\fcZ{5,1}\equiv -\fcZ{5,1},\fcZR{1,5}$														\\\cdashline{1-2}
		$\wspZVd{6}{2} + \wspPZV{6}$     	&$\fcZ{4,1,1}\equiv -\fcZ{4,1,1},\pm\fcZR{1,4,1},\fcZR{1,1,4}$										\\
								     	&$\fcZ{3,2,1}\equiv -\fcZ{2,3,1},\fcZR{1,3,2}$													\\
								     	&$\fcZ{3,1,2}\equiv -\fcZ{3,1,2},\fcZ{2,1,3}$													\\
								     	&$\fcZ{2,3,1}\equiv -\fcZ{3,2,1},\fcZR{1,2,3}$													\\
								     	&$\fcZ{2,1,3}\equiv \fcZ{3,1,2},-\fcZ{2,1,3}$													\\\cdashline{1-2}
		$\wspZVd{6}{3} + \wspPZV{6}$     	&$\fcZ{3,1,1,1}\equiv -\fcZ{3,1,1,1},\pm\fcZR{1,3,1,1},\pm\fcZR{1,1,3,1},\fcZR{1,1,1,3}$					\\
								   	&$\fcZ{2,2,1,1}\equiv -\fcZ{2,2,1,1},\pm\fcZR{1,2,2,1},\fcZR{1,1,2,2},$								\\
								   	&$\fcZ{2,1,2,1}\equiv -\fcZ{2,1,2,1},\pm\fcZ{2,1,1,2},\fcZR{1,2,1,2}$									\\
								   	&$\fcZ{2,1,1,2}\equiv \pm\fcZ{2,1,2,1},\pm\fcZR{1,2,1,2}$											\\\cdashline{1-2}
		$\wspZVd{6}{4} + \wspPZV{6}$     	&$\fcZ{2,1,1,1,1}\equiv -\fcZ{2,1,1,1,1},\pm\fcZR{1,2,1,1,1},\pm\fcZR{1,1,2,1,1},\pm\fcZR{1,1,1,2,1},\fcZR{1,1,1,1,2}$	\\\hline
	\end{tabular} 
	\caption{Examples of \refThm{1_Thm2}}\label{3_PL_TableExamples}
}\end{center}
\end{table}

We list examples of \refThm{1_Thm2} for MZVs with weight $l$ smaller than $7$ in \refTab{3_PL_TableExamples}.
As an application of the examples,
	we will give a property of MZVs of weight $6$ such that
	\envHLineCmPt[3_Ex1_EqMain]{\equiv}
	{
		\fcZ{l_1,\ldots,l_n}
	}
	{
		0
		\qquad\modd \wspZVd{6}{n-1} + \wspPZV{6}
	}
	where $n\in\Set{2,3,4,5}$, $l_1+\cdots+l_n=6$ and  $l_1\geq2$.
For this,
	we may prove  
	\envHLineCFCmPt{\equiv}
	{
		0
	}
	{
		\fcZ{3,3}
		\hspace{127pt}\qquad \modd \wspZVd{6}{1} + \wspPZV{6}
	}
	{
		0
	}
	{
		\fcZ{3,2,1}
	\lnP{\equiv}
		\fcZ{2,3,1}
	\lnP{\equiv}
		\fcZ{2,2,2}
		\qquad \modd \wspZVd{6}{2} + \wspPZV{6}
	}
	because the other cases are obvious from $\fcZ{l_1,\ldots,l_n}\equiv-\fcZ{l_1,\ldots,l_n}$ (see \refTab{3_PL_TableExamples}).
By the harmonic relations obtained by \refEq{2_PL_DefHarProdMain},
	we have
	\envM{
		2\fcZ{3,3} 
	}{
		\fcZ{3}\fcZ{3} - \fcZ{6}
	}
	and
	\envMCm{
		3\fcZ{2,2,2}
	}
	{
		\fcZ{2}\fcZ{2,2} - \fcZ{4,2} - \fcZ{2,4} 
	}
	which verify 
	\envMPt{\equiv}{
		\fcZ{3,3} 
	}{
		0	\ \modd \wspZVd{6}{1} + \wspPZV{6}
	}
	and
	\envMCmPt{\equiv}{
		\fcZ{2,2,2}
	}
	{
		0	\ \modd \wspZVd{6}{2} + \wspPZV{6}
	}
	respectively.
By the duality theorem \refEq{3_Lem2Proof_EqDuality},
	we have $\fcZ{2,3,1}=\fcZ{3,1,2}$,
	where the space giving the congruence relation $\equiv$ is $\wspZVd{6}{2} + \wspPZV{6}$.
Since
	$\fcZ{2,3,1}\equiv-\fcZ{3,2,1}$ and $\fcZ{3,1,2}\equiv0$,
	we obtain
	$\fcZ{3,2,1}\equiv\fcZ{2,3,1}\equiv0$,
	which completes the proof of \refEq{3_Ex1_EqMain}.





\end{document}